\documentclass[a4paper,10pt]{scrartcl}
\pdfoutput=1
\usepackage{ifpdf}
\usepackage[english]{babel}
\usepackage[utf8]{inputenc}
\usepackage{enumerate}
\usepackage{geometry}
\usepackage{amsfonts,amsmath,amssymb,amsopn}
\usepackage{amsthm}
\usepackage{mathtools}
\usepackage{stmaryrd}
\usepackage{nicefrac}
\usepackage{exscale}

\usepackage[numbers,square]{natbib}
\usepackage{url} 
\usepackage[colorlinks=true,pdfpagelabels,unicode]{hyperref}
\hypersetup{linkcolor=blue, urlcolor=blue, citecolor=blue} 

\allowdisplaybreaks 

\theoremstyle{plain}

\newtheoremstyle{theo}
	{3pt} 
	{3pt} 
	{\itshape} 
	{} 
		{\bfseries} 
	{\\} 
	{ } 
	{\thmname{#1}\thmnumber{ #2.}\thmnote{ - #3}} 
\theoremstyle{theo}

\newtheorem{definition}{Definition}[section]
\newtheorem{lemma}[definition]{Lemma}
\newtheorem{theorem}[definition]{Theorem}
\newtheorem{corollary}[definition]{Corollary}
\newtheorem{proposition}[definition]{Proposition}

\newenvironment{bew}{\begin{proof}[\bfseries Proof:]}{\end{proof}}

\newtheoremstyle{remark}
	{3pt} 
	{3pt} 
	{} 
	{} 
		{\bfseries} 
	{} 
	{ } 
	{\thmname{#1}\thmnumber{ #2.}\thmnote{ - #3}} 
\theoremstyle{remark}
\newtheorem{remark}[definition]{Remark}

\DeclareMathOperator{\bomega}{\overline{\Omega}}
\DeclareMathOperator{\romega}{\partial\Omega}

\DeclareMathOperator{\intd}{d\!}

\DeclareMathOperator{\wto}{\rightharpoonup}
\DeclareMathOperator{\wsto}{\stackrel{\star}{\wto}}

\newcommand{\epsi}{\varepsilon}

\usepackage[usenames,dvipsnames,svgnames,table]{xcolor}

\newcommand{\Ts}{T_{\star}}

\newcommand{\Td}{T_{\diamond}}
\newcommand{\Tme}{T_{max,\;\!\epsi}}
\newcommand{\GNI}{Gagliardo--Nirenberg inequality}

\newcommand{\into}[1]{\int_0^{#1}\!}
\newcommand{\intot}{\into{t}}
\newcommand{\intoT}{\into{T}}

\newcommand{\intomega}{\int_{\Omega}\!} 

\newcommand{\intoTomega}{\intoT\!\intomega}
\newcommand{\intotomega}{\intot\!\intomega}
\newcommand{\intinfomega}{\int_0^\infty\!\!\intomega}

\newcommand{\Lo}[1][1]{L^{#1}(\Omega)} 
\newcommand{\W}[1][1,2]{W^{#1}(\Omega)}

\newcommand{\LSp}[2]{L^{#1\;\!}\!\left(#2\right)} 

\newcommand{\LSploc}[2]{L_{loc}^{#1}\!\left(#2\right)} 

\newcommand{\CSp}[2]{C^{#1}\!\left(#2\right)}

\newcommand{\CSpnl}[2]{C^{#1}\!\,(#2)} 

\newcommand{\R}{\mathbb{R}}
\newcommand{\N}{\mathbb{N}}

\newcommand{\uep}{u_\epsi}
\newcommand{\vep}{v_\epsi}
\newcommand{\wep}{w_\epsi}
\newcommand{\ubs}{{u_*}}
\newcommand{\vbs}{{v_*}}
\newcommand{\ws}{w_*}
\newcommand{\rs}{r_{*}}
\newcommand{\rhs}{{r}_{**}}

\author{Tobias Black\thanks{Institut f\"ur Mathematik, Universit\"at Paderborn, Warburger Str. 100, 33098 Paderborn, Germany; email: \mbox{tblack@math.upb.de}}}
\title{Global generalized solutions to a forager-exploiter model with superlinear degradation and their eventual regularity properties}
\date{}

\begin{document}
\maketitle
\begin{abstract}
\noindent
{\textbf{Abstract:}
\noindent In this article we consider a cascaded taxis model for two proliferating and degrading species which thrive on the same nutrient but orient their movement according to different schemes. In particular, we assume the first group, the foragers, to orient their movement directly along an increasing gradient of the food density, while the second group, the exploiters, instead track higher densities of the forager group. Specifically, we will investigate an initial boundary-value problem for a prototypical forager-exploiter model of the form
\begin{align*}
\left\{
\begin{array}{r@{\,}c@{\,}l@{\quad }l@{\quad}l@{\quad}l@{\,}c}
u_{t}&=&\Delta u-\nabla\cdot\big(u\nabla w\big)+f(u),\ &x\in\Omega,& t>0,\\
v_{t}&=&\Delta v-\nabla\cdot\big(v\nabla u\big)+g(v),\ &x\in\Omega,& t>0,\\
w_{t}&=&\Delta w-(u+v)w-\mu w+r(x,t),\ &x\in\Omega,& t>0,
\end{array}\right.
\end{align*}
in a smoothly bounded domain $\Omega\subset\mathbb{R}^2$, where $\mu\geq 0$, $r\in\ C^1\!\left(\overline{\Omega}\times[0,\infty)\right)\cap L^\infty\!\left(\Omega\times(0,\infty)\right)$ is nonnegative and the functions $f,g\in C^1\!\left([0,\infty)\right)$ are assumed to satisfy $f(0)\geq0$, $g(0)\geq0$ as well as $$-k_f s^{\alpha} -l_f\leq f(s)\leq -K_f s^\alpha+L_f\ \text{ and }\ -k_g s^\beta -l_g\leq g(s)\leq -K_g s^\beta+L_g\quad\text{for }s\geq0,$$ respectively, with constants $\alpha,\beta>1$, $k_f,K_f,k_g,K_g>0$ and $l_f,L_f,l_g,L_g\geq0$ and $\alpha,\beta>1$.  

Assuming that $\alpha>1+\sqrt{2}$, $\min\{\alpha,\beta\}>\frac{\alpha+1}{\alpha-1}$ and that $r$ satisfies certain structural conditions we establish the global solvability of this system with respect to a suitable generalized solution concept and then, for the more restrictive case of $\alpha,\beta>1+\sqrt{2}$ and $\mu>0$, investigate an eventual regularity effect driven by the decay of the nutrient density $w$.
}\\[0.1cm]


\end{abstract}


\newpage
\section{Introduction}\label{sec1:intro}
{\noindent
Interaction between groups of the same or even different species following distinct taxis schemes to adapt their movement is known to support rich dynamical behavior and the emergence of spatial patterns. Examples can not only be found in biological contexts of collective migration, such as swarming and flocking (\cite{eftimie07,tonerFlocksHerdsSchools1998,bellomo12}), but also in other fields like economy (\cite{furioli17}) and criminology (\cite{short08}). A particular instance, where only two different groups are involved in the process, can be witnessed within Alaska's bird population with the macroscopic formation of shearwater flocks which orient themselves towards kittiwake foragers in order to find feeding grounds with a sufficient food source (\cite{hoffmann81}). Observations like this are the essential motivation behind the typical forager-exploiter nutrient models, in which one forager group actively searches for the food source, while the second exploiter group only indirectly finds the nutrients by tracking the actively searching foragers.

To describe a corresponding minimal model the authors of \cite{tania12} proposed a cascaded taxis system essentially of the form
\begin{align}\label{for-expl-no-source}
\left\{
\begin{array}{r@{\,}c@{\,}l@{\quad }l@{\quad}l}
u_{t}&=&\Delta u-\nabla\cdot\big(u\nabla w\big),\\
v_{t}&=&\Delta v-\nabla\cdot\big(v\nabla u\big),\\
w_{t}&=&\Delta w-(u+v)w-\mu w+r.
\end{array}\right.
\end{align}
Herein, $u=u(x,t)$ denotes the time-evolution of the density of the forager population, $v=v(x,t)$ the density of the exploiter population, $w=w(x,t)$ the density of their nutrients. The degradation rate $\mu\geq0$ is assumed to be constant and the external resupply of nutrients $r=r(x,t)$ is supposed to be a sufficiently regular nonnegative function. From an application point of view, however, it seems even more appropriate to incorporate the possibilities of degradation and proliferation for the forager and exploiter populations into the model, where, in particular, on larger time-scales for the motivating example of shearwater flocks and kittiwake foragers the death of individuals at high population densities should not be neglected. Following the approaches of related settings, this consideration would usually be integrated into the system by introducing logistic source terms of the forms $+au-bu^2$ and $+av-bv^2$ with $a,b>0$ to the first and second equations, respectively, where the prominent role of the exponent $2$ is further underlined by its appearance in the extensively studied Fisher-KPP equation (\cite{Fisher37,Kolmo91}). In this work, nevertheless, we are going to consider more general systems of the form
\begin{align}\label{forager-0}
\left\{
\begin{array}{r@{\,}c@{\,}l@{\quad }l@{\quad}l}
u_{t}&=&\Delta u-\nabla\cdot\big(u\nabla w\big)+f(u),\\
v_{t}&=&\Delta v-\nabla\cdot\big(v\nabla u\big)+g(v),\\
w_{t}&=&\Delta w-(u+v)w-\mu w+r,
\end{array}\right.
\end{align}
where we will assume that the prescribed kinetic source terms $f,g$ only satisfy certain growth restrictions (see \eqref{fgprop} below), which ensure degradation for large population densities and thereby include spontaneous death effects.

From a mathematical point of view a cascaded taxis-mechanism of the form in \eqref{for-expl-no-source} and \eqref{forager-0} is quite challenging as even a single taxis term already allows for a considerable amount of mathematical questions, e.g. illustrated by the extensive studies of the acclaimed Keller--Segel chemotaxis model. (See the surveys \cite{BBWT15,HP09,LanWin_JDMV_19} for an overview of related models and results.) In order to better understand the analytical difficulty in this framework let us first take a look at the closely related setting with only one homogeneous group, where a archetypal model, often referred to as prey-taxis system (\cite{leePatternFormationPreytaxis2009,WuShiWu16}) or chemotaxis-consumption system (\cite{Tao-consumption_JMAA11,TaoWin12_evsmooth,LanWan-Highdim-chemo-cons_17}) in the literature, can be given in the form
\begin{align}\label{prey-taxis}
\left\{
\begin{array}{r@{\,}c@{\,}l@{\quad }l@{\quad}l}
u_{t}&=&\Delta u-\chi\nabla\cdot\big(u\nabla w\big),\\
w_{t}&=&\Delta w-uw-\mu w+r,
\end{array}\right.
\end{align}
where the parameter of chemotactic strength $\chi>0$ and the degradation rate $\mu\geq0$ are constant and $r=r(x,t)\geq0$ is sufficiently regular. Again $u=u(x,t)$ and $w=w(x,t)$ denote the density of the population and their nutrients, respectively. The prototypical case of $\mu=r=0$ has been studied rigorously and global existence and boundedness of solutions is only known under smallness conditions on the initial data (\cite{Tao-consumption_JMAA11}), or in a two-dimensional setting (\cite{Win-glob-large-data_CPDE12,liGlobalExistenceUniform2015}), where also their stabilization towards an homogeneous equilibrium has been established for all reasonably regular initial data by utilizing the energy structure present in the system (\cite{win-stab2d-ArchRatMechAna12}). Even in settings where additional source terms are introduced to \eqref{prey-taxis}, suitable adaptations of the energy method can provide quite extensive insight in the corresponding systems as e.g. witnessed by the results on global existence, boundedness and global attractors for suitably small $\chi>0$ in \cite{WuShiWu16}.

In the cascaded system \eqref{for-expl-no-source}, however, the favorable energy structure of \eqref{prey-taxis} is missing due to the sequential attraction of the two populations and, accordingly, the analytical knowledge of the system is still quite sparse and often only covers low dimensions or generalized solution concepts. Particularly, in the one-dimensional setting the authors of \cite{taoLargeTimeBehavior2019} established global classical solutions for suitably regular initial data and, additionally, proved an exponential stabilization result under the assumption that one of the initial masses of $u$ and $v$ is sufficiently small. In higher dimensions, however, only generalized solution concepts could be considered, with \cite{Win-globgensocialint-M3AS18} providing results on global generalized solutions for $n\geq1$ whenever the initial datum for $w$ satisfies a certain smallness condition. These solutions, moreover, approach spatially homogeneous profiles in the large time limit in some topology if the supply function $r$ decays sufficiently fast. A recent study, where a volume-filling effect was included by replacing $u\nabla w$ with $u(1-u)\nabla w$ and $v\nabla u$ with $v(1-v)\nabla u$ in the first and second equation of \eqref{for-expl-no-source}, respectively, established global bounded classical solutions in general dimension under the assumption that the initial data for $u$ and $v$ are uniformly bounded by the density threshold $1$ (\cite{liuGlobalExistenceBoundedness2019}). 

On the other hand, in single species chemotaxis-growth systems it could be observed that including a sufficiently degrading source term can have a distinct regularizing effect on the solution components. Notably, the chemotaxis systems with logistic growth terms
\begin{align}\label{KScons}
\left\{
\begin{array}{r@{\,}c@{\,}l@{\quad }l@{\quad}l}
u_{t}&=&\Delta u-\chi\nabla\cdot\big(u\nabla w\big)+au-bu^2,\\
w_{t}&=&\Delta w-uw,
\end{array}\right.
\end{align}
with $a\geq0$ and $b>0$ admit global classical solutions in general dimensions if $b$ is large enough compared to the initial datum of $w$ and $\chi$ (\cite{LanWan-Highdim-chemo-cons_17}). In three-dimensional chemotaxis-Navier--Stokes systems the logarithmic growth term, moreover, fuels an eventual smoothing effect even when $b>0$ is arbitrary small (\cite{Lan-Longterm_M3AS16}), further highlighting the regularizing influence emitted by kinetic source terms. This is even more dominant in the classical Keller--Segel setting, where the attracting substance is produced by the population instead of consumed, i.e. replacing the second equation by $w_t=\Delta w-w+u$ in \eqref{KScons}. This setting with $f(u)\equiv0$ allows for blow-up solutions to exist in dimension $n=2$ if the initial mass of $u$ is large enough (\cite{HVblow97,HW01,mizoguchi_winkler_13}) and for arbitrary mass if $n\geq3$ (\cite{Win13pure}), while all solutions are global and bounded whenever $n=2$ and $f(u)=au-bu^2$ with $a\geq0$ and $b>0$ or whenever $n\geq3$ and $b>0$ is sufficiently large (\cite{win10ctax}). Accordingly, source functions with a dominating death term at high population densities appear to be favorable for the global existence of solutions and we have high hopes that in the physically meaningful setting of \eqref{forager-0} in a two-dimensional domain we can obtain distinctively better results beyond the quite restrictive global existence results for \eqref{for-expl-no-source} mentioned above.

\noindent\textbf{Main results.} For the remainder of the work we are going to investigate the initial-boundary value problem
\begin{align}\label{forager}
\left\{
\begin{array}{r@{\,}c@{\,}l@{\quad }l@{\quad}l@{\quad}l@{\,}c}
u_{t}&=&\Delta u-\nabla\cdot\big(u\nabla w\big)+f(u),\ &x\in\Omega,& t>0,\\
v_{t}&=&\Delta v-\nabla\cdot\big(v\nabla u\big)+g(v),\ &x\in\Omega,& t>0,\\
w_{t}&=&\Delta w-(u+v)w-\mu w+r(x,t),\ &x\in\Omega,& t>0,\\
\nabla u\cdot\nu&=&0,\quad \nabla v\cdot\nu=0,\quad \nabla w\cdot\nu=0,\ &x\in\romega,&t>0,\\
u(x,0)&=&u_0(x),\ v(x,0)=v_0(x),\ w(x,0)=w_0(x),\ &x\in\Omega,&
\end{array}\right.
\end{align}
where $\nu$ denotes the outward normal vector field on $\romega$ and $f,g\in\CSp{1}{[0,\infty)}$ satisfy 
\begin{equation}\label{fgprop}
\begin{aligned}
-k_f s^{\alpha} -l_f\leq f(s)\leq -K_f s^\alpha+L_f&\ \text{ and }\ -k_g s^\beta -l_g\leq g(s)\leq -K_g s^\beta+L_g,\quad\text{for }s\geq 0\\
&\text{as well as}\quad f(0)\geq0,\quad g(0)\geq0,\\
\end{aligned}
\end{equation}
with constants $K_f,K_g,k_f,k_g>0$, $L_f,L_g,l_f,l_g\geq 0$ and $\alpha,\beta> 1$. As for the other ingredients of \eqref{forager} we merely require that $\mu$ is a nonnegative constant, that $r\in\CSp{1}{\bomega\times[0,\infty)}\cap\LSp{\infty}{\Omega\times(0,\infty)}$ is nonnegative with 
\begin{align}\label{rprop}
\rs:=\sup_{s\in(0,\infty)}\|r(\cdot,s)\|_{\Lo[\infty]}<\infty
\end{align}
and that the initial data $(u_0,v_0,w_0)$ satisfy the conditions
\begin{align}\label{IR}
\left\{\begin{array}{r@{\ }l}
 u_0\in\W[1,\infty]\ &\text{is nonnegative with }\ u_0\not\equiv0,\\
 v_0\in\W[1,\infty]&\text{is nonnegative with }\ v_0\not\equiv0,\\
 w_0\in\W[2,\infty]&\text{is nonnegative with }\ \nabla w_0\cdot\nu=0\text{ on }\romega.
\end{array}\right. 
\end{align}

In this setting our first result is concerned with the global existence of generalized solutions for kinetic functions $f$ and $g$ satisfying \eqref{fgprop} with suitably large degradation rates $\alpha,\beta>1$ and reads as follows. 

\begin{theorem}\label{theo:1} 
Let $\Omega\subset\R^2$ be a bounded domain with smooth boundary. Suppose that $\mu\geq0$, that the functions $f$ and $g$ fulfill \eqref{fgprop} with $K_f,K_g,k_f,k_g>0$, $L_f,L_g,l_f,l_g\geq0$, $\alpha>1+\sqrt{2}$ and $\beta>1$ satisfying
\begin{align}\label{eq:min-cond}
\min\{\alpha,\beta\}>\frac{\alpha+1}{\alpha-1}
\end{align}
and that $r\in\CSp{1}{\bomega\times[0,\infty)}\cap\LSp{\infty}{\Omega\times(0,\infty)}$ is nonnegative and fulfills \eqref{rprop}. Then, for any $u_0,v_0$ and $w_0$ complying with \eqref{IR}, the system \eqref{forager} admits at least one global generalized solution $(u,v,w)$ in the sense of Definition \ref{def:sol}.
\end{theorem}

The prominent degradation exponent $2$, which is commonly found in studies of systems with only a single taxis term, is not covered by the result above and hence the question whether in the cascaded taxis system logistic growth terms are strong enough to ensure global solutions in quite mild solution concepts has to be left for further research. Nevertheless, physically meaningful source terms, e.g. functions describing an \emph{Allee-effect}, are still contained in the theorem above. We remark the following.  

\begin{remark}\label{rem:1}\mbox{}
\begin{enumerate}
\item[i)] We note that if $\alpha\leq\beta$, condition \eqref{eq:min-cond} is always satisfied due to $\alpha>1+\sqrt{2}$,  whereas for large $\alpha>\beta$ the fraction $\frac{\alpha+1}{\alpha-1}$ tends to $1$ and hence if $\alpha$ is suitably large, the degradation function $g$ may also be of a subquadratic but superlinear form.
\item[ii)] Kinetic source terms like $f(s)\equiv s(1-s)(s-2)$, which are prototypical choices when describing a population evolution obeying certain \emph{Allee-effects}, satisfy the conditions of Theorem \ref{theo:1}. 
\end{enumerate}
\end{remark}

In our second result we investigate a slightly more restrictive setting, where we assume that in addition to all the requirements of Theorem \ref{theo:1} we also have $\beta>1+\sqrt{2}$, $\mu>0$ and that $r\in\CSp{1}{\bomega\times[0,\infty)}\cap\LSp{\infty}{\Omega\times(0,\infty)}\cap\LSp{1}{(0,\infty);\Lo[\infty]}$ satisfies
\begin{align}\label{rprop2}
\rhs:=\int_0^\infty\!\|r(\cdot,t)\|_{\Lo[\infty]}\intd t<\infty.
\end{align}

These conditions allow for an eventual regularization process to occur, which entails that after some (possibly large) waiting time $\Td>0$ the destabilizing effect of the cascaded taxis terms is mitigated by the sufficiently strong degradation of the forager and exploiter groups and the suitably fast decaying resupply of nutrients implied by \eqref{rprop2}, so that the solution to \eqref{forager} actually becomes a classical solution after some time.

\begin{theorem}\label{theo:2} 
Let $\Omega\subset\R^2$ be a bounded domain with smooth boundary. Suppose that $\mu>0$, that the functions $f$ and $g$ fulfill \eqref{fgprop} with $K_f,K_g,k_f,k_g>0$, $L_f,L_g,l_f,l_g\geq0$, $\alpha>1+\sqrt{2}$ and $\beta>1+\sqrt{2}$ and that $r\in\CSp{1}{\bomega\times[0,\infty)}\cap\LSp{\infty}{\Omega\times(0,\infty)}\cap\LSp{1}{(0,\infty);\Lo[\infty]}$ is nonnegative and fulfills \eqref{rprop} and \eqref{rprop2}. Then there exists $\Td>0$ such that the global generalized solution of \eqref{forager} from Theorem \ref{theo:1} has the properties
\begin{align*}
u\in\CSp{2,1}{\bomega\times[\Td,\infty)},\qquad v\in\CSp{2,1}{\bomega\times[\Td,\infty)},\quad\text{and}\quad w\in\CSp{2,1}{\bomega\times[\Td,\infty)},
\end{align*}
and such that $(u,v,w)$ solves \eqref{forager} classically in $\Omega\times(\Td,\infty)$.
\end{theorem}

\noindent\textbf{Outline of the approach.} After a brief introduction of the generalized solution concept we are going to consider (Definition \ref{def:sol}), we will propose a family of regularized initial boundary-value problems in Section \ref{sec3:prelim} and, starting with a time-local existence result, provide uniform $L^1$ bounds for both $\uep$ and $\vep$ and a uniform $L^\infty$ bound on $\wep$, which, when combined with semigroup arguments, entail the time-global existence of the approximating system for any fixed $\epsi\in(0,1)$ (Lemma \ref{lem:globex}). Section \ref{sec4:further-estimates} is concerned with the derivation of essential a priori estimates, which constitute the main part of the proof of Theorem \ref{theo:1}. In particular, we will draw on the well-known maximal Sobolev regularity estimates for the Neumann heat semigroup employed to the third equation (Lemma \ref{lem:spat-temp-est-W2r-wep}) to obtain sufficient control on the terms appearing during the testing procedure for $\uep$ (Lemma \ref{lem:st-bounds-uep-q}). The second quantity will be treated by considering the time-evolution of $\intomega\ln(\vep+1)$ (Lemma \ref{lem:st-bounds-vep}), which offers very mild, yet still sufficient, information on the spatial gradient of $\vep$. In Section \ref{sec5:equiintegrability-estimates} we prepare additional a priori estimates necessary for the limiting procedure undertaken in Proposition \ref{prop:conv} of Section \ref{sec6:limit}. In Section \ref{sec7:solutionprop} we will then conclude the proof of Theorem \ref{theo:1} by verifying that the obtained limit object from Section \ref{sec6:limit} is in fact a global generalized solution.

The proof of Theorem \ref{theo:2} will be in the focus of Sections \ref{sec8:ev-cond-sob-reg}-\ref{sec11:ev-smooth}, with Section \ref{sec8:ev-cond-sob-reg} being dedicated to conditional regularity estimates which provide the basis for the bootstrapping procedure undertaken in Section \ref{sec11:ev-smooth}. Sections \ref{sec9:w-decay} and \ref{sec10:ev-lp-uep} ensure that the conditions required by the conditional estimates of Lemma \ref{lem:ev-max-sob} are indeed fulfilled. Specifically, we make use of \eqref{rprop2} and the fact that $\mu>0$ to prove the decay of $\wep$ in the $L^\infty$ norm in Lemma \ref{lem:dec-wep-Linfty} and then investigate $\frac{\intd}{\intd t}\intomega\frac{\uep^q}{(2\delta-\wep)^\theta}$ for suitably small $\delta>0$ and $\theta>0$ to establish an eventual $L^q$ bound for $\uep$ in Lemma \ref{lem:ev-lp-bound-uep}. Section \ref{sec11:ev-smooth} finalizes the proof of Theorem \ref{theo:2} by the aforementioned bootstrapping procedure based on Lemma \ref{lem:ev-max-sob}.
}

\setcounter{equation}{0} 
\section{Global solutions in a generalized solution concept}\label{sec2:soldef}
Since both of our theorems are concerned with global generalized solutions, let us first specify what exactly constitutes a solution in this notion. The main difference to the standard concept of weak solvability consists of only requiring a integral inequality corresponding to the variational interpretation of the second equation of \eqref{forager}. This in turn allows us to prescribe milder regularity assumptions for the second solution component, which certainly is the hardest to obtain for the system in question. Solution concepts of this form have been successfully employed in many chemotaxis contexts, see e.g. \cite{win15_chemorot,LanWin2017}.

\begin{definition}\label{def:sol}
A triple $(u,v,w)$ of functions
\begin{align*}
u&\in\LSploc{2}{[0,\infty);\W[1,2]},\quad\ 
v\in\LSploc{1}{\bomega\times[0,\infty)},\quad\ 
w\in\LSploc{\infty}{\bomega\times[0,\infty)}\cap\LSploc{2}{[0,\infty);\W[1,2]}
\end{align*}
such that
\begin{align*}
f(u)\quad&\text{and}\quad\frac{g(v)}{v+1}\quad\text{belong to }\LSploc{1}{\bomega\times[0,\infty)}\quad
\text{and that}\quad \nabla\ln(v+1)\in\LSploc{2}{\bomega\times[0,\infty);\R^2},
\end{align*}
and such that $u\geq0$, $v\geq0$ and $w\geq0$ in $\bomega\times[0,\infty)$ will be called a global generalized solution of \eqref{forager} if the inequality
\begin{align}\label{eq:def-ineq-v}
-\intinfomega&\ln(v+1)\psi_t-\intomega\ln(v_0+1)\psi(\cdot,0)\nonumber\\&\geq \intinfomega|\nabla\ln(v+1)|^2\psi-\intinfomega\nabla\ln(v+1)\cdot\nabla\psi-\intinfomega\frac{v}{v+1}\big(\nabla\ln(v+1)\cdot\nabla u\big)\psi\\&\hspace*{1.6cm}+\intinfomega\frac{v}{v+1}\nabla u\cdot\nabla\psi+\intinfomega\frac{g(v)}{v+1}\psi\nonumber
\end{align}
holds for every nonnegative $\psi\in C^\infty_0(\bomega\times[0,\infty))$, if the identities
\begin{align}\label{eq:def-eq-u}
-\intinfomega u\varphi_t-\intomega u_0\varphi(\cdot,0)=-\intinfomega \nabla u\cdot\nabla\varphi+\intinfomega u\nabla w\cdot\nabla\varphi+\intinfomega f(u)\varphi
\end{align}
and
\begin{align}\label{eq:def-eq-w}
-\intinfomega w\varphi_t-\intomega w_0\varphi(\cdot,0)=-\intinfomega\nabla w\cdot\nabla\varphi-\intinfomega(u+v)w\varphi-\mu\intinfomega w\varphi+\intinfomega r\varphi
\end{align}
are fulfilled for every $\varphi\in C_0^\infty(\bomega\times[0,\infty))$ and if finally the inequality
\begin{align}\label{eq:def-mass-ineq}
\intomega v(\cdot,t)\leq \intomega v_0+\intotomega g(v)
\end{align}
holds for a.e. $t\in(0,\infty)$.
\end{definition}

We note that an estimate of the form \eqref{eq:def-mass-ineq} is a necessary requirement, if we want to have a concept of generalized solvability, which is consistent with the concept of classical solvability in the following sense.

\begin{remark}\label{rem:2}
Whenever a global generalized solution $(u,v,w)$ satisfies $u,v,w\in\CSp{0}{\bomega\times[0,\infty)}\cap\CSp{2,1}{\bomega\times(0,\infty)}$, it can easily be checked that the global generalized solution is also a global classical solution of \eqref{forager} in $\Omega\times(0,\infty)$. (See \cite[Lemma 2.1]{win15_chemorot} and \cite[Lemma 2.5]{LanWin2017} for the essential arguments involved.)
\end{remark}

\setcounter{equation}{0} 
\section{Global solutions to approximating systems and basic properties}\label{sec3:prelim}

Global generalized solutions in the sense of the Definition \ref{def:sol} above will be obtained as limit of global classical solutions to the following family of regularized systems
\begin{align}\label{approxprob}
\left\{
\begin{array}{r@{\,}c@{\,}l@{\quad }l@{\quad}l@{\quad}l@{\,}c}
u_{\epsi t}&=&\Delta \uep-\nabla\cdot\big(\uep\nabla \wep\big)+f(\uep),\ &x\in\Omega,& t>0,\\
v_{\epsi t}&=&\Delta \vep-\nabla\cdot\big(\vep\nabla \uep\big)+g(\vep),\ &x\in\Omega,& t>0,\\
w_{\epsi t}&=&\Delta \wep-\frac{(\uep+\vep)\wep}{1+\epsi(\uep+\vep)\wep}-\mu \wep+r(x,t),\ &x\in\Omega,& t>0,\\
\nabla \uep\cdot\nu&=&0,\quad \nabla \vep\cdot\nu=0,\quad \nabla \wep\cdot\nu=0,\ &x\in\romega,&t>0,\\
\uep(x,0)&=&u_0(x),\ \vep(x,0)=v_0(x),\ \wep(x,0)=w_0(x),\ &x\in\Omega,&
\end{array}\right.
\end{align}
where in particular the fact that $\frac{\xi}{1+\epsi\xi}\leq\frac{1}{\epsi}$ for all $\epsi\in(0,1)$ and $\xi\geq0$ enables us to obtain quite strong regularity estimates on the third solution component for each fixed $\epsi\in(0,1)$, providing an important tool when trying to show that for each $\epsi\in(0,1)$ the solution to \eqref{approxprob} is global in time. As we still need a few additional ingredients, however, let us begin our analysis of this family of approximating systems with merely establishing time-local existence of classical solutions to \eqref{approxprob}.

\begin{lemma}
\label{lem:locex}
Let $\Omega\subset\R^2$ be a bounded domain with smooth boundary. Suppose that the functions $f$ and $g$ fulfill \eqref{fgprop} with $K_f,K_g,k_f,k_g>0$, $L_f,L_g,l_f,l_g\geq0$ and some $\alpha,\beta> 1$. Assume $\mu\geq0$ and that $r\in\CSp{1}{\bomega\times[0,\infty)}\cap\LSp{\infty}{\Omega\times(0,\infty)}$ is nonnegative and satisfies \eqref{rprop}. Then, for any $u_0,v_0$ and $w_0$ complying with \eqref{IR} and any $\epsi\in(0,1)$ there exist $\Tme>0$ and functions
\begin{align}\label{eq:reg-approx-sol}
\begin{cases}
\uep\in \bigcap_{q>n}\CSp{0}{[0,\Tme);\W[1,q]}\cap\CSp{2,1}{\bomega\times[0,\Tme)},\\
\vep\in \bigcap_{q>n}\CSp{0}{[0,\Tme);\W[1,q]}\cap\CSp{2,1}{\bomega\times[0,\Tme)},\\
\wep\in \bigcap_{q>n}\CSp{0}{[0,\Tme);\W[1,q]}\cap\CSp{2,1}{\bomega\times[0,\Tme)},
\end{cases}
\end{align}
which solve \eqref{forager} in $\Omega\times(0,\Tme)$ in the classical sense and satisfy $\uep\geq0$, $\vep\geq 0$ and $\wep\geq 0$ in $\bomega\times(0,\infty)$. Moreover, either $\Tme=\infty$ or
\begin{align}\label{eq:extensibility}
\limsup_{t\nearrow\Tme}\Big(\|\uep(\cdot,t)\|_{\W[1,p]}+\|\vep(\cdot,t)\|_{\W[1,p]}+\|\wep(\cdot,t)\|_{\W[1,p]}\Big)=\infty
\end{align}
for all $p>2$.
\end{lemma}

\begin{bew}
The time-local existence of a classical solution $(\uep,\vep,\wep)$ satisfying \eqref{eq:reg-approx-sol} is a consequence of Amann's general theory on quasilinear parabolic boundary value problems (cf. \cite[Theorem 14.5 and Theorem 14.6]{Amann93-Nonhomogeneoues-linear-and-quasilinear}) with the extensibility criterion \eqref{eq:extensibility} entailed by \cite[Theorem 15.5]{Amann93-Nonhomogeneoues-linear-and-quasilinear}. The nonnegativity of the solution components is then ensured by employing the parabolic comparison principle to the subsolution 0.
\end{bew}

For the remainder of this work we will always assume that all conditions of Lemma \ref{lem:locex} are satisfied and for $\epsi\in(0,1)$ we will always denote by $(\uep,\vep,\wep)$ the solution provided by Lemma \ref{lem:locex} and with $\Tme\in(0,\infty]$ its corresponding maximal existence time. 

Relying on our standing assumption \eqref{fgprop}, we can integrate the first and second equations of \eqref{approxprob} to obtain an $L^1$ bound uniform in $\epsi\in(0,1)$ for both $\uep$ and $\vep$.
\begin{lemma}
\label{lem:L1-est}
For any $\epsi\in(0,1)$ and all $t\in(0,\Tme)$ the solution $(\uep,\vep,\wep)$ of \eqref{approxprob} obtained in Lemma \ref{lem:locex} satisfies
\begin{align}\label{eq:L1-est}
\intomega\uep(\cdot,t)\leq\max\bigg\{\intomega u_0,\ |\Omega|\Big(\frac{L_f}{K_f}\Big)^\frac{1}{\alpha}\bigg\}=:\ubs\quad\text{and}\quad \intomega\vep(\cdot,t)\leq\max\bigg\{\intomega v_0,\ |\Omega|\Big(\frac{L_g}{K_g}\Big)^\frac{1}{\beta}\bigg\}=:\vbs.
\end{align}
\end{lemma}

\begin{bew}
Making use of the first equation and the prescribed boundary conditions of \eqref{approxprob}, we can draw on integration by parts and \eqref{fgprop} to find that
\begin{align*}
\frac{\intd}{\intd t}\intomega \uep=\intomega f(\uep)\leq -K_f\intomega \uep^\alpha+L_f|\Omega|
\end{align*}
holds on $(0,\Tme)$. Due to $\alpha>1$ an application of Hölder's inequality entails that $y_\epsi(t):=\intomega \uep(\cdot,t)$ satisfies the differential inequality
\begin{align*}
y_\epsi'(t)+\frac{K_f}{|\Omega|^{\alpha-1}} y_\epsi^\alpha(t)\leq L_f|\Omega|\quad\text{for all }t\in(0,\Tme),
\end{align*}
An ODE comparison argument with $\ubs$ then yields the asserted $L^1$ bound for $\uep$. Repeating similar arguments for the second equation of \eqref{approxprob} we also conclude the desired bound for $\intomega \vep$.
\end{bew}
Emulating arguments applied to the system without kinetic source functions, we can make use of the same supersolution as in \cite[Lemma 3.1]{Win-globgensocialint-M3AS18} to easily conclude the uniform $L^\infty$ boundedness of $\wep$ by a straightforward comparison argument.

\begin{lemma}
\label{lem:Linfty-est-w}
For any $\epsi\in(0,1)$ and all $t\in(0,\Tme)$ the solution $(\uep,\vep,\wep)$ of \eqref{approxprob} satisfies
\begin{align*}
\|\wep(\cdot,t)\|_{\Lo[\infty]}\leq \|w_0\|_{\Lo[\infty]}+\frac{\rs}{\mu}=:\ws.
\end{align*}
\end{lemma}

\begin{bew}
For $x\in\bomega$ and $t\in(0,\infty)$ we let
\begin{align*}
\overline{w}(x,t):=\|w_0\|_{\Lo[\infty]}e^{-\mu t}+\int_0^t e^{-\mu(t-s)}\|r(\cdot,s)\|_{\Lo[\infty]}\intd s.
\end{align*}
Then, we observe that $\overline{w}$ is spatially homogeneous for each $t\in(0,\Tme)$ with $\overline{w}\leq \ws$ on $\bomega\times(0,\Tme)$ and satisfies
\begin{align*}
\overline{w}_t-\Delta\overline{w}+\frac{(\uep+\vep)\wep}{1+\epsi(\uep+\vep)\wep}+\mu\overline{w}-r(x,t)\geq\overline{w}_t+\mu\overline{w}-\|r(\cdot,t)\|_{\Lo[\infty]}=0\quad\text{on }\Omega\times(0,\Tme),
\end{align*}
due to the nonnegativity of $\uep,\vep$ and $\wep$. Hence, an application of the comparison principle entails that $\overline{w}\geq \wep$ on $\Omega\times(0,\Tme)$ and, in particular, $\|\wep(\cdot,t)\|_{\Lo[\infty]}\leq \ws$ for all $t\in(0,\Tme)$.
\end{bew}

The $L^\infty$ information for $\wep$ at hand we can now extract sufficient information to verify that for any fixed $\epsi\in(0,1)$ the solution is global in time.

\begin{lemma}
\label{lem:globex}
For any $\epsi\in(0,1)$ the solution $(\uep,\vep,\wep)$ of \eqref{approxprob} obtained in Lemma \ref{lem:locex} is global, i.e. $\Tme=\infty$.
\end{lemma}

\begin{bew}
We augment the arguments presented in \cite[Lemma 2.1]{Win-globgensocialint-M3AS18} to cover the present setting. We fix a suitably large $p>2$ and assume that $\Tme<\infty$. Due to $\frac{\xi}{1+\epsi\xi}\leq\frac{1}{\epsi}$ for all $\xi\geq0$ and $r\leq\rs$ in $\Omega\times(0,\Tme)$ we can on the one hand utilize well-known maximal Sobolev regularity theory (see e.g. \cite{HieberPruss-CPDE92,GigSohr91}), Lemma \ref{lem:Linfty-est-w} and \eqref{IR} to obtain $C_1>0$ such that
\begin{align*}
&\int_0^{\Tme}\|w_{\epsi t}\|_{\Lo[p]}^p+\int_0^{\Tme}\|\wep\|_{\W[2,p]}^p\\\leq\ &C_1\|w_0\|_{\W[2,p]}^p+C_1\int_0^{\Tme}\Big(\Big\|\frac{(\uep+\vep)\wep}{1+\epsi(\uep+\vep)\wep}\Big\|_{\Lo[p]}^p+\|\mu\wep\|_{\Lo[p]}^p+\|r\|_{\Lo[p]}^p\Big)\\
\leq\ &C_1\big(\|w_0\|_{\W[2,p]}^p+\epsi^{-p}|\Omega|+\mu^p\ws^p|\Omega|+\rs^p|\Omega|\big)\Tme,
\end{align*}
and on the other hand employ semigroup techniques (e.g. \cite[Lemma 1.3]{win10jde}) to obtain $C_2>0$ such that also $\|\nabla\wep\|_{\Lo[\infty]}\leq C_2$ holds on $(0,\Tme)$. This information at hand, we perform testing procedures on the first equation of \eqref{approxprob} and rely on the fact that $\alpha>1$ to obtain $C_3>0$ such that that $\|\uep\|_{\Lo[4]}\leq C_3$ on $(0,\Tme)$, which, by utilizing the semigroup estimates for the Neumann heat semigroup once more, can also be refined into a bound for $\|\uep\|_{\Lo[\infty]}$ on $(0,\Tme)$. In light of the properties assumed for $f$, the boundedness of $\uep$ also readily entails the existence of $C_4>0$ such that $\|f(\uep)\|_{\Lo[\infty]}\leq C_4$ on $(0,\Tme)$. Hence, we may interpret $\uep$ as a weak solution of the equation
\begin{align*}
\tilde{u}_t=\Delta\tilde{u}+a(x,t)\cdot\nabla\tilde{u}+b(x,t)\tilde{u}+c(x,t)\quad\text{on }\Omega\times(0,\Tme),
\end{align*}
with $a(x,t):=-\nabla\wep\in\LSp{\infty}{\Omega\times(0,\Tme);\R^2}$, $b(x,t)=-\Delta\wep\in\LSp{p}{\Omega\times(0,\Tme)}$ and $c(x,t):=f(\uep)\in\LSp{\infty}{\Omega\times(0,\Tme)}$, where we can now draw on standard theory to conclude that also $\uep\in\LSp{\infty}{(0,\Tme);\W[1,\infty]}\cap\LSp{p}{(0,\Tme);\W[2,p]}$. With this we may repeat the steps employed to $\uep$ again for the second equation of \eqref{approxprob} to obtain $\vep\in\LSp{\infty}{(0,\Tme);\W[1,\infty]}\cap\LSp{p}{(0,\Tme);\W[2,p]}$, which finally contradicts the extensibility criterion \eqref{eq:extensibility} and hence our assumption $\Tme<\infty$ has to be false.
\end{bew}

\setcounter{equation}{0} 
\section{Essential uniform estimates}\label{sec4:further-estimates}
This section is devoted to the establishment of estimates which build the essential groundwork for later refinement in Section \ref{sec5:equiintegrability-estimates} and the limit procedure in Section \ref{sec6:limit}. This section will also feature the main reasons behind the constraints imposed on the parameters $\alpha$ and $\beta$ in Theorem \ref{theo:1}.

We start by noting that an integration of the first and second equations of \eqref{approxprob} also entails some additional information we did not mention in Lemma \ref{lem:L1-est}. In particular, we also obtain an inequality for $\vep$ corresponding to \eqref{eq:def-mass-ineq} of the generalized solution concept.
\begin{lemma}
\label{lem:spat-temp-alpha-beta}
Let $\ubs$ and $\vbs$ be provided by \eqref{eq:L1-est}. For any $\epsi\in(0,1)$ the solution $(\uep,\vep,\wep)$ of \eqref{approxprob} fulfills
\begin{align}\label{eq:spat-temp-alpha-beta}
\int_t^{t+1}\!\!\intomega\uep^\alpha\leq \frac{L_f|\Omega|}{K_f}+\frac{\ubs}{K_f}\quad\text{and}\quad\int_t^{t+1}\!\!\intomega\vep^\beta\leq\frac{L_g|\Omega|}{K_g}+\frac{\vbs}{K_g}
\end{align}
for all $t>0$. Moreover, for any $\epsi\in(0,1)$ and all $t>0$ the solution satisfies
\begin{align}\label{eq:spat-temp-minus-plus}
\intomega\vep(\cdot,t)=\intomega v_0+\intotomega g(\vep).
\end{align}
\end{lemma}

\begin{bew}
Again by integrating the first and second equations of \eqref{approxprob} over $\Omega$ we first find that
\begin{align}\label{eq:spat-temp-alpha-beta-1}
\frac{\intd}{\intd t}\intomega\uep=\intomega f(\uep)\quad\text{and}\quad\frac{\intd}{\intd t}\intomega\vep=\intomega g(\vep)\quad\text{hold on }(0,\infty).
\end{align}
Plugging in \eqref{fgprop} we obtain that
\begin{align*}
\frac{\intd}{\intd t}\intomega\uep+K_f\intomega\uep^\alpha\leq L_f|\Omega|\quad\text{and}\quad\frac{\intd}{\intd t}\intomega\vep+K_g\intomega\vep^\beta=L_g|\Omega|
\end{align*}
are valid on $(0,\infty)$. Since $\uep$ and $\vep$ are nonnegative, we conclude upon integration from $t$ to $t+1$ that
\begin{align*}
&K_f\int_t^{t+1}\!\!\intomega\uep^\alpha\leq \intomega\uep(\cdot,t)+L_f|\Omega|\quad\text{and}\quad K_g\int_t^{t+1}\!\!\intomega\vep^\beta\leq \intomega\vep(\cdot,t)+L_g|\Omega|
\end{align*}
hold for any $t>0$, which, in light of Lemma \ref{lem:L1-est}, proves \eqref{eq:spat-temp-alpha-beta}.
Going back to \eqref{eq:spat-temp-alpha-beta-1}, we observe that integrating with respect to time immediately entails \eqref{eq:spat-temp-minus-plus}, concluding the proof.
\end{bew}
In order to prepare for the core argument of our analysis in the next lemma we will require the following temporal cut-off functions. 

\begin{definition}
\label{def:cutoff}
Let $\xi_0\in\CSp{\infty}{\R}$ be an arbitrary monotonically increasing function satisfying
\begin{align*}
0\leq\xi_0\leq 1\ \text{ on }\ \R,\qquad\xi_0\equiv0\ \text{ on }\ (-\infty,0],\quad\text{and}\quad \xi_0\equiv 1\ \text{ on }\ (1,\infty).
\end{align*}
For $T>0$ we then define
\begin{align*}
\xi_T(t):=\xi_0(t-T),\quad t\in\R.
\end{align*}
\end{definition}

The essential driving force of our analysis will be obtained by means of a well-known maximal Sobolev regularity estimate for the Neumann heat semigroup (\cite{HieberPruss-CPDE92,GigSohr91}) applied to the third equation of \eqref{approxprob}. In our setting, relying on the result of Lemma \ref{lem:spat-temp-alpha-beta}, we are hence able to obtain time-space bounds for the second derivatives of $\wep$ in $L^q$ spaces with $1<q<\min\{\alpha,\beta\}$.

\begin{lemma}
\label{lem:spat-temp-est-W2r-wep}
Let $\rho:=\min\{\alpha,\beta\}>1$. There exists $C>0$ such that for any $\epsi\in(0,1)$ and all $t\geq0$ the solution $(\uep,\vep,\wep)$ of \eqref{approxprob} satisfies
\begin{align*}
\int_t^{t+1}\!\|\wep\|_{\W[2,\rho]}^\rho\leq C.
\end{align*}
\end{lemma}

\begin{bew}
For arbitrary fixed $\tau>0$ we let $\xi:=\xi_\tau$ denote the cut-off function introduced in Definition \ref{def:cutoff} and note that $\xi\wep$ solves 
\begin{align*}
\big(\xi\wep\big)_t&=\Delta(\xi\wep)-\frac{(\uep+\vep)(\xi\wep)}{1+\epsi(\uep+\vep)\wep}-\mu(\xi\wep)+\xi r(x,t)+\xi'\wep\quad\text{in }\Omega\times(\tau,\infty)\\
&\text{with }\big(\xi\wep\big)(\cdot,\tau)=0\quad\text{in }\Omega\quad\text{and}\quad \nabla(\xi\wep)\cdot\nu=0\quad\text{on }\romega\times(\tau,\infty).
\end{align*}
Hence, according to the maximal Sobolev regularity estimate for the Neumann heat semigroup (\cite{HieberPruss-CPDE92,GigSohr91}) 
 for any $q>1$ there exists some $C_1>0$ such that for any $\epsi\in(0,1)$ the inequality
\begin{align*}
\int_\tau^{\tau+2}\!\|(\xi\wep)_t\|_{\Lo[\rho]}^\rho&+\int_\tau^{\tau+2}\!\|\xi\wep\|_{\W[2,\rho]}^\rho\\
&\leq C_1\bigg(0+ \int_\tau^{\tau+2}\!\Big\|-\frac{(\uep+\vep)\xi\wep}{1+\epsi(\uep+\vep)\wep}-\mu\xi\wep+ \xi r+\xi'\wep\Big\|_{\Lo[\rho]}^\rho\bigg)
\end{align*}
holds. Since $\frac{1}{1+\epsi(\uep+\vep)\wep}\leq 1$ on $\Omega\times(0,\infty)$, we can draw on the nonnegativity of $r,\uep,\vep$ and $\wep$ on $\Omega\times(0,\infty)$, the fact that $\xi\leq 1$ on $\R$, Lemma \ref{lem:Linfty-est-w} and \eqref{rprop} to obtain $C_2>0$ such that
\begin{align*}
\int_\tau^{\tau+2}\!\|\xi\wep\|_{\W[2,\rho]}^\rho&\leq 2C_1|\Omega|\big(\mu^\rho\ws^\rho+\rs^\rho+\|\xi'\|_{\LSp{\infty}{\R}}^\rho\ws^\rho\big)+C_1\ws^q\int_\tau^{\tau+2}\!\!\intomega(\uep+\vep)^\rho\\
&\leq C_2+C_2\int_\tau^{\tau+2}\!\!\intomega\uep^\rho+C_2\int_\tau^{\tau+2}\!\!\intomega\vep^\rho\quad\text{for all }\epsi\in(0,1).
\end{align*}
Here, we rely on \eqref{eq:spat-temp-alpha-beta} from Lemma \ref{lem:spat-temp-alpha-beta} to find that
\begin{align*}
\int_\tau^{\tau+2}\!\|\xi\wep\|_{\W[2,\rho]}^\rho&\leq C_2+2C_2\left(\frac{L_f|\Omega|}{K_f}+\frac{\ubs}{K_f}+\frac{L_g|\Omega|}{K_g}+\frac{\vbs}{K_g}\right),
\end{align*}
which, due to the arbitrary choice of $\tau>0$ and the fact that $\xi\equiv 1$ on $(\tau+1,\tau+2)$, implies the existence of $C_3>0$ such that
\begin{align*}
\int_t^{t+1}\!\|\wep\|_{\W[2,\rho]}^\rho\leq C_3
\end{align*} 
holds for each $\epsi\in(0,1)$ and all $t>1$. We are left with treating $t\in[0,1]$. Arguing along similar lines, this time investigating \eqref{approxprob} directly without the cut-off and additionally relying on \eqref{IR} to treat the first term, we find $C_4>0$ such that for any $\epsi\in(0,1)$
\begin{align*}
\int_0^{2}\!\|\wep\|_{\W[2,\rho]}^\rho&\leq C_1\|w_0\|_{\W[2,\rho]}^\rho+2C_1|\Omega|\big(\mu^\rho\ws^\rho+\rs^\rho\big)+ C_1\ws^\rho\int_0^{2}\!\!\intomega(\uep+\vep)^\rho\\
&\leq C_4+C_2\int_0^{2}\!\!\intomega\big(\uep^\rho+\vep^\rho\big)\leq C_4+2C_2\left(\frac{L_f|\Omega|}{K_f}+\frac{\ubs}{K_f}+\frac{L_g|\Omega|}{K_g}+\frac{\vbs}{K_g}\right),
\end{align*}
completing the proof upon obvious choice for $C>0$. 
\end{bew}

In the next step we will see that the information on the Laplacian of $\wep$, contained in Lemma \ref{lem:spat-temp-est-W2r-wep}, is the key ingredient when making use of testing procedures to derive $L^q$ bounds for $\uep$. The following lemma is also the crucial point where the conditions on $\alpha$ and $\min\{\alpha,\beta\}$ featured in Theorem \ref{theo:1} are first required.

\begin{lemma}
\label{lem:testing-uep}
Assume $\alpha>1+\sqrt{2}$ and $\beta>1$ and suppose that $\rho:=\min\{\alpha,\beta\}$ satisfies $\rho>\frac{\alpha+1}{\alpha-1}$. Then for any $q\in\big[2,(\alpha-1)(\rho-1)\big)$ there exists $C>0$ such that for all $\epsi\in(0,1)$ the solution $(\uep,\vep,\wep)$ of \eqref{approxprob} satisfies
\begin{align*}
\frac{1}{q}\frac{\intd}{\intd t}\intomega \uep^q+(q-1)\intomega\uep^{q-2}|\nabla\uep|^2+\frac{K_f}{2}\intomega\uep^{q+\alpha-1}\leq C\intomega|\Delta\wep|^\rho+C
\end{align*}
on $(0,\infty)$.
\end{lemma}

\begin{bew}
First note that, due to the assumption that $\rho>\frac{\alpha+1}{\alpha-1}$, the interval $I:=\big[2,(\alpha-1)(\rho-1)\big)$ is not empty. Now, given any $q\in I$ we test the first equation of \eqref{approxprob} against $\uep^{q-1}$ and integrate by parts, which implies that
\begin{align*}
\frac{1}{q}\frac{\intd}{\intd t}\intomega\uep^q=-(q-1)\intomega \uep^{q-2}|\nabla \uep|^2+(q-1)\intomega\uep^{q-1}(\nabla\uep\cdot\nabla\wep)+\intomega\uep^{q-1} f(\uep)
\end{align*}
holds on $(0,\infty)$, due to the prescribed boundary conditions. In light of the pointwise identity $\uep^{q-1}\nabla\uep=\frac{1}{q}\nabla(\uep^q)$ we may integrate by parts once more in the second term on the right to obtain
\begin{align*}
\frac{1}{q}\frac{\intd}{\intd t}\intomega\uep^q=-(q-1)\intomega\uep^{q-2}|\nabla\uep|^2-\frac{q-1}{q}\intomega\uep^q\Delta\wep+\intomega\uep^{q-1} f(\uep)
\end{align*}
is valid on $(0,\infty)$. Here, we make use of the nonnegativity of $\uep$, \eqref{fgprop} and Young's inequality to find $C_1>0$ such that for any $\epsi\in(0,1)$
\begin{align*}
\frac{1}{q}\frac{\intd}{\intd t}\intomega\uep^q\leq-(q-1)\intomega\uep^{q-2}|\nabla\uep|^2+C_1\intomega|\Delta\wep|^\rho+\intomega\uep^{\frac{q\rho}{\rho-1}}-K_f\intomega\uep^{q+\alpha-1}+L_f\intomega\uep^{q-1}
\end{align*}
on $(0,\infty)$. Since $\alpha$ is nonnegative and $q<(\alpha-1)(\rho-1)$ implies that also $\frac{q\rho}{\rho-1}<q-1+\alpha$ we may employ Young's inequality once again to the third and last terms on the right and conclude that there is some $C_2>0$ such that for any $\epsi\in(0,1)$ the inequality
\begin{align*}
\frac{1}{q}\frac{\intd}{\intd t}\intomega\uep^q+(q-1)\intomega\uep^{q-2}|\nabla\uep|^2+\frac{K_f}{2}\intomega \uep^{q+\alpha-1}\leq C_1\intomega|\Delta\wep|^\rho+C_2
\end{align*}
is satisfied on $(0,\infty)$.
\end{bew}

In order to turn this differential inequality into some boundedness information we require the following supplementary lemma. This lemma is taken from \cite[Lemma 3.4]{LL18-ClassSolLogCTSingSens}, where also additional details on the proof can be found.

\begin{lemma}
\label{lem:diffineq-lemma}
For some $T\in(0,\infty]$ let $y\in\CSp{1}{(0,T)}\cap\CSp{0}{[0,T)}$, $h\in\CSp{0}{[0,T)}$, $h\geq0$, $C>0$, $a>0$ satisfy
\begin{align*}
y'(t)+ay(t)\leq h(t),\qquad\int_{(t-1)_+}^t h(s)\intd s\leq C
\end{align*}
for all $t\in(0,T)$. Then $y\leq y(0)+\frac{C}{1-e^{-a}}$ throughout $(0,T)$.
\end{lemma}

A combination of Lemmas \ref{lem:testing-uep} and \ref{lem:diffineq-lemma} with the boundedness of $\int_t^{t+1}\!\!\intomega|\Delta\wep|^\rho$ provided by Lemma \ref{lem:spat-temp-est-W2r-wep} now entails uniform (spatio-temporal) bounds for $\uep$ in some $L^q$ spaces, where the admissible values for $q$ are already sufficiently large for our later purposes.

\begin{lemma}
\label{lem:st-bounds-uep-q}
Assume $\alpha>1+\sqrt{2}$ and $\beta>1$ and suppose that $\rho:=\min\{\alpha,\beta\}$ satisfies $\rho>\frac{\alpha+1}{\alpha-1}$. Then for any $q\in\big[2,(\alpha-1)(\rho-1)\big)$ there exists $C>0$ such that for all $\epsi\in(0,1)$ the solution $(\uep,\vep,\wep)$ of \eqref{approxprob} fulfills
\begin{align*}
\intomega \uep^q(\cdot,t)+\int_t^{t+1}\!\!\intomega\uep^{q-2}|\nabla\uep|^2+\int_t^{t+1}\!\!\intomega\uep^{q+\alpha-1}\leq C
\end{align*}
for all $t>0$.
\end{lemma}

\begin{bew}
From Lemma \ref{lem:testing-uep} we obtain $C_1>0$ such that
\begin{align}\label{eq:st-bounds-uep-eq1}
\frac{1}{q}\frac{\intd}{\intd t}\intomega \uep^q+(q-1)\intomega\uep^{q-2}|\nabla\uep|^2+\frac{K_f}{2}\intomega\uep^{q+\alpha-1}\leq C_1\intomega|\Delta\wep|^\rho+C_1
\end{align}
is valid on $(0,\infty)$ for any $\epsi\in(0,1)$. To derive an inequality of the form featured in Lemma \ref{lem:diffineq-lemma} from this, we note that by Lemma \ref{lem:spat-temp-est-W2r-wep} there is $C_2>0$ satisfying
\begin{align}\label{eq:st-bounds-uep-eq0}
\int_t^{t+1}\!\!\intomega|\Delta\wep|^\rho\leq C_2\quad\text{for all }\epsi\in(0,1)\text{ and }t>0,
\end{align}
and that $\alpha>1+\sqrt{2}$ and Young's inequality entail the existence of $C_3=C_3(q,\alpha,K_f,\Omega)>0$ such that $\frac{1}{q}\intomega\uep^q(\cdot,t)\leq\frac{K_f}{2}\intomega\uep^{q+\alpha-1}(\cdot,t)+C_3$ for all $t>0$. Hence, setting $y_\epsi(t):=\frac{1}{q}\intomega\uep^q(\cdot,t)$, $t>0$, we have
\begin{align*}
y_\epsi'(t)+ y_\epsi(t)\leq  C_1\intomega|\Delta\wep|^\rho+ C_1+ C_3\quad\text{for all }\epsi\in(0,1)\text{ and }t>0,
\end{align*}
where we used the nonnegativity of $\uep$ to drop the term containing the derivative. Now we can apply Lemma \ref{lem:diffineq-lemma} in conjunction with \eqref{eq:st-bounds-uep-eq0} to obtain $C_4>0$ such that
\begin{align*}
\intomega\uep^q(\cdot,t)\leq C_4\quad\text{for all }\epsi\in(0,1)\text{ and }t>0.
\end{align*}
Finally, we see from \eqref{eq:st-bounds-uep-eq1} upon integration over $(t,t+1)$ that
\begin{align*}
\frac{1}{q}\intomega\uep^q(\cdot,t+1)+(q-1)\int_t^{t+1}\!\!\intomega\uep^{q-2}|\nabla\uep|^2+\frac{K_f}{2}\int_t^{t+1}\!\!\intomega\uep^{q+\alpha-1}\leq C_4+ C_1C_2+C_1\quad\text{for all }t>0,
\end{align*}
concluding the proof.
\end{bew}

An evident consequence of the lemma above are the following bounds.
\begin{corollary}
\label{cor:st-bounds-uep-l2}
Assume $\alpha>1+\sqrt{2}$ and $\beta>1$ and suppose that $\rho:=\min\{\alpha,\beta\}$ satisfies $\rho>\frac{\alpha+1}{\alpha-1}$. Then there exists $C>0$ such that for all $\epsi\in(0,1)$ the solution $(\uep,\vep,\wep)$ of \eqref{approxprob} fulfills
\begin{align*}
\intomega \uep^2(\cdot,t)+\int_t^{t+1}\!\!\intomega|\nabla\uep|^2+\int_t^{t+1}\!\!\intomega\uep^{\alpha+1}\leq C
\end{align*}
for all $t>0$.
\end{corollary}

\begin{bew}
This follows immediately from Lemma \ref{lem:st-bounds-uep-q} when applied for $q=2$.
\end{bew}

To establish the uniform bound for $\uep$ in $\Lo[2]$ we made essential use of the knowledge that $\Delta\wep$ is uniformly bounded in $\LSp{\rho}{\Omega\times(t,t+1)}$. Since up to this point the only regularity information on any derivative of $\uep$ is a mere bound on $\nabla\uep$ in $\LSp{2}{\Omega\times(t,t+1)}$, it is therefore not surprising that the treatment of $\vep$ is exceedingly difficult and we cannot expect standard testing measures to improve our knowledge on $\vep$ beyond the result of Lemma \ref{lem:spat-temp-alpha-beta}. Hence, we turn to the investigation of the rather mild expression $\ln(\vep+1)$, which will at least provide minimal information on the spatial derivative of $\vep$ in the sense of the following lemma.

\begin{lemma}
\label{lem:st-bounds-vep}
Assume $\alpha>1+\sqrt{2}$ and $\beta>1$ and suppose that $\rho:=\min\{\alpha,\beta\}$ satisfies $\rho>\frac{\alpha+1}{\alpha-1}$. Then for any $T>0$ there exists $C(T)>0$ such that for all $\epsi\in(0,1)$ the solution $(\uep,\vep,\wep)$ of \eqref{approxprob} satisfies
\begin{align*}
\int_0^{T}\!\!\intomega\frac{|\nabla\vep|^2}{(\vep+1)^2}\leq C(T).
\end{align*}
\end{lemma}

\begin{bew}
Straightforward calculations using the second equation of \eqref{approxprob} and integration by parts show that
\begin{align*}
-\frac{\intd}{\intd t}\intomega\ln(\vep+1)&=-\intomega\frac{\Delta\vep-\nabla\cdot(\vep\nabla\uep)+g(\vep)}{\vep+1}\\
&=-\intomega\frac{|\nabla\vep|^2}{(\vep+1)^2}+\intomega\frac{\vep}{(\vep+1)^2}(\nabla\vep\cdot\nabla\uep)-\intomega\frac{g(\vep)}{\vep+1}
\end{align*}
is valid for all $t>0$. Since $\vep$ is nonnegative, an application of Young's inequality and \eqref{fgprop} thereby entails that
\begin{align*}
-\frac{\intd}{\intd t}\intomega\ln(\vep+1)+\frac{1}{2}\intomega\frac{|\nabla\vep|^2}{(1+\vep)^2}&\leq-\intomega g(\vep)+ \frac{1}{2}\intomega\frac{\vep^2}{(1+\vep)^2}|\nabla\uep|^2\\\
&\leq k_g\intomega\vep^\beta+l_g|\Omega|+\frac{1}{2}\intomega|\nabla\uep|^2\quad\text{for all }t>0.
\end{align*}
Given any $T>0$, we can now integrate this inequality with respect to time and immediately obtain that
\begin{align}\label{eq:st-bounds-vep-eq1}
\frac{1}{2}\int_0^T\!\!\intomega\frac{|\nabla\vep|^2}{(\vep+1)^2}\leq \intomega\ln\big(\vep(\cdot,T)+1\big)+l_g|\Omega|T+k_g\int_0^T\!\!\intomega\vep^\beta+\int_0^T\!\!\intomega|\nabla\uep|^2.
\end{align}
To further estimate the right hand side, we recall that Lemma \ref{lem:spat-temp-alpha-beta} and Corollary \ref{cor:st-bounds-uep-l2} provide $C_1(T)>0$ such that $k_g\int_0^T\!\!\intomega\vep^\beta +\int_0^T\!\!\intomega|\nabla\uep|^2\leq C_1(T)$ and that Lemma \ref{lem:L1-est}, when combined with the the evident estimate $0\leq\ln(\xi+1)\leq \xi$ for all $\xi\geq0$, entails $\intomega\ln(\vep(\cdot,T)+1)\leq\intomega\vep(\cdot,T)\leq \vbs$. Plugging these bounds into \eqref{eq:st-bounds-vep-eq1} finally implies
\begin{align*}
\frac{1}{2}\int_0^T\!\!\frac{|\nabla\vep|^2}{(\vep+1)^2}\leq \vbs+l_g|\Omega|T+C_1(T),
\end{align*}
completing the proof.
\end{bew}

\setcounter{equation}{0} 
\section{Refined spatio-temporal estimates and the regularity of time derivatives -- Preparing convergence in suitable \texorpdfstring{$L^p$}{Lp} spaces}\label{sec5:equiintegrability-estimates}
The goal of this section is the preparation of uniform bounds in $L^p$ spaces for suitably chosen $p>1$. In the succeding sections we will combine these bounds with Vitali's theorem to attain the strong convergence properties necessary for passing to limit in the integrals appearing in integral identities corresponding to \eqref{eq:def-ineq-v}, \eqref{eq:def-eq-u} and \eqref{eq:def-eq-w} of the solution concept (see Section \ref{sec7:solutionprop}). Let us begin with a straightforward implication of Lemma \ref{lem:spat-temp-est-W2r-wep}.

\begin{lemma}
\label{lem:vitali-nab-w}
Let $\rho:=\min\{\alpha,\beta\}>1$. There is $C>0$ such that for all $t>0$ and all $\epsi\in(0,1)$ the solution $(\uep,\vep,\wep)$ of \eqref{approxprob} satisfies
\begin{align*}
\int_t^{t+1}\!\!\intomega|\nabla \wep|^{2\rho}\leq C.
\end{align*}
\end{lemma}

\begin{bew}
The \GNI\ asserts the existence of $C_1>0$ such that
\begin{align*}
\|\phi\|_{\W[1,2\rho]}\leq C_1\|\phi\|_{\W[2,\rho]}^a\|\phi\|_{\Lo[\infty]}^{1-a}\quad\text{for all }\phi\in\W[2,\rho],
\end{align*}
where $a=\frac{\frac{1}{2}-\frac{1}{2\rho}}{1-\frac{1}{\rho}}=\frac{1}{2}$. Making use of Lemma \ref{lem:Linfty-est-w} we thus obtain that $\wep$ satisfies
\begin{align*}
\int_t^{t+1}\!\|\wep\|_{\W[1,2\rho]}^{2\rho}\leq C_1\int_t^{t+1}\!\|\wep\|_{\W[2,\rho]}^\rho\|\wep\|_{\Lo[\infty]}^\rho\leq C_1\ws^\rho\int_t^{t+1}\!\|\wep\|_{\W[2,\rho]}^\rho
\end{align*}
for all $t>0$ and all $\epsi\in(0,1)$. In light of Lemma \ref{lem:spat-temp-est-W2r-wep} this entails the existence of $C_2>0$ such that
\begin{align*}
\int_t^{t+1}\!\|\wep\|_{\W[1,2\rho]}^{2\rho}\leq C_2\quad\text{for all }t>0\text{ and }\epsi\in(0,1),
\end{align*}
from which we immediately conclude the asserted bound.
\end{bew}

Another direct implication of the bounds prepared at the beginning of Section \ref{sec4:further-estimates} and the $L^\infty$ estimate for $\wep$ is the following lemma.
\begin{lemma}
\label{lem:st-bound-l1-uepvepwep}
Let $\rho:=\min\{\alpha,\beta\}>1$. For all $s\in(1,\rho)$ and $T>0$ there is $C(T)>0$ such that for all $\epsi\in(0,1)$ the solution $(\uep,\vep,\wep)$ of \eqref{approxprob} fulfills
\begin{align*}
\int_0^{T}\!\!\intomega\Big(\frac{(\uep+\vep)\wep}{1+\epsi(\uep+\vep)\wep}\Big)^s\leq C(T).
\end{align*}
\end{lemma}

\begin{bew}
We fix any $s\in(1,\rho)$. Since $\frac{1}{(1+\epsi(\uep+\vep)\wep)^s}\leq 1$ on $\Omega\times(0,\infty)$, we can make use of Young's inequality and Lemma \ref{lem:Linfty-est-w} to estimate
\begin{align*}
\int_0^T\!\!\intomega\Big(\frac{(\uep+\vep)\wep}{1+\epsi(\uep+\vep)\wep}\Big)^s\leq\int_0^T\!\!\intomega(\uep+\vep)^\rho+\int_0^T\!\!\intomega\wep^{\frac{s\rho}{\rho-s}}\leq 2^\rho\int_0^T\!\!\intomega \uep^\rho+2^\rho\int_0^T\!\!\intomega\vep^\rho+\ws^\frac{s\rho}{\rho-s}|\Omega| T.
\end{align*}
The boundedness of the remaining terms is a direct consequence of Lemma \ref{lem:spat-temp-alpha-beta} due to $\rho=\min\{\alpha,\beta\}$.
\end{bew}
In the next lemma we will take a look at bounds targeting $f(\uep)$ and $g(\vep)$. Here we note that, in particular, the uniform bound for $f(\uep)$ in $\LSp{\frac{(\alpha+1)s}{\alpha}}{\Omega\times(0,T)}$ for some $s>1$ is one essential part in order to later refine the convergence result for $\nabla\uep$ to strong convergence in $\LSp{2}{\Omega\times(0,T)}$ (cf. Section \ref{sec7:solutionprop}), which is also the main reason we had to treat the permissible exponents in Lemma \ref{lem:testing-uep} so carefully.

\begin{lemma}
\label{lem:equiint-f-g}
Assume $\alpha>1+\sqrt{2}$ and $\beta>1$ and suppose that $\rho:=\min\{\alpha,\beta\}$ satisfies $\rho>\frac{\alpha+1}{\alpha-1}$. Then there exists $s>1$ such that for any $T>0$ there is $C(T)>0$ such that for all $\epsi\in(0,1)$ the solution $(\uep,\vep,\wep)$ of \eqref{approxprob} fulfills
\begin{align*}
\int_0^{T}\!\!\intomega\big|f(\uep)\big|^{\frac{\alpha+1}{\alpha}s}+\int_0^{T}\!\!\intomega\Big|\frac{g(\vep)}{\vep+1}\Big|^{s}\leq C(T).
\end{align*}
\end{lemma}

\begin{bew}
We set $s_0:=\min\big\{\frac{\alpha-1}{\alpha+1}\rho,\frac{\beta}{\beta-1}\big\}$ and note that the conditions $\rho>\frac{\alpha+1}{\alpha-1}$ and $\beta>1$ clearly imply that $s_0>1$. Hence, we can pick some $s\in(1,s_0)$ and note that since $(\alpha+1)s<(\alpha-1)\rho$ and $(\beta-1)s<\beta$ we can rely on Lemma \ref{lem:st-bounds-uep-q}, employed for $q=(\alpha+1)s-(\alpha-1)<(\alpha-1)(\rho-1)$, and Lemma \ref{lem:spat-temp-alpha-beta} to obtain $C_1(T)>0$ such that
\begin{align}\label{eq:equiint-fg-eq1}
\int_0^T\!\!\intomega\uep^{(\alpha+1)s}+\int_0^T\!\!\intomega\vep^{(\beta-1)s}\leq C_1(T)\quad\text{for all }\epsi\in(0,1).
\end{align}
 Now, making use of \eqref{fgprop} we find $C_2>0$ such that
\begin{align*}
\int_0^T\!\!\intomega\big|f(\uep)\big|^{\frac{\alpha+1}{\alpha}s}+\int_0^T\!\!\intomega\Big|\frac{g(\vep)}{\vep+1}\Big|^s&\leq C_2\int_0^T\!\!\intomega \uep^{(\alpha+1)s}+C_2\int_0^T\!\!\intomega\frac{\vep^{\beta s}}{(\vep+1)^s}+C_2 T
\end{align*}
holds for any $T>0$ and all $\epsi\in(0,1)$. Relying on \eqref{eq:equiint-fg-eq1} and the fact that $\xi^{\beta s}(\xi+1)^{-s}\leq\xi^{(\beta-1)s}$ for all $\xi\geq0$ we find that
\begin{align*}
\int_0^T\!\!\intomega\big|f(\uep)\big|^{\frac{\alpha+1}{\alpha}s}+\int_0^T\!\!\intomega\Big|\frac{g(\vep)}{\vep+1}\Big|^s&\leq C_2C_1(T)+C_2C_1(T)+C_2 T
\end{align*}
is valid for any $T>0$ and all $\epsi\in(0,1)$.
\end{bew}

Making use of Corollary \ref{cor:st-bounds-uep-l2} and Lemma \ref{lem:vitali-nab-w} we can also readily attain a $\epsi$-uniform bound on $\uep\nabla\wep$ in $\LSp{s}{\Omega\times(0,T)}$ for some $s>2$. This result will also be the second essential part for establishing the strong convergence of $\nabla\uep$ in Section \ref{sec7:solutionprop}.
\begin{lemma}
\label{lem:vitali-u-nabw}
Assume $\alpha> 1+\sqrt{2}$ and $\beta>1$ and suppose that $\rho:=\min\{\alpha,\beta\}$ satisfies $\rho>\frac{\alpha+1}{\alpha-1}$. Then there exists some $s>2$ such that for any $T>0$ there is $C(T)>0$ such that for all $\epsi\in(0,1)$ the solutions $(\uep,\vep,\wep)$ of \eqref{approxprob} satisfies
\begin{align*}
\int_0^T\!\!\intomega|\uep\nabla \wep|^s\leq C(T).
\end{align*}
\end{lemma}

\begin{bew}
We note that due to $\rho>\frac{\alpha+1}{\alpha-1}$ the interval $I:=\big(2,\frac{2\rho(\alpha+1)}{\alpha+1+2\rho}\big)$ is not empty and that for $s\in I$ we have $s<\alpha+1$ and $\frac{(\alpha+1)s}{\alpha+1-s}< 2\rho$. Accordingly, for any fixed $s\in I$ straightforward applications of Young's inequality entail the existence of $C_1>0$ such that
\begin{align*}
\int_0^T\!\!\intomega\uep^s|\nabla\wep|^s\leq\int_0^T\!\!\intomega\uep^{\alpha+1}+\int_0^T\!\!\intomega|\nabla\wep|^{\frac{(\alpha+1)s}{\alpha+1-s}}\leq \int_0^T\!\!\intomega\uep^{\alpha+1}+\int_0^T\!\!\intomega|\nabla\wep|^{2\rho}+C_1 T
\end{align*}
holds for any $T>0$ and all $\epsi\in(0,1)$. In light of the previously prepared bounds of Corollary \ref{cor:st-bounds-uep-l2} and Lemma \ref{lem:vitali-nab-w} we can immediately conclude the assertion.
\end{bew}

Before we can employ an Aubin--Lions type argument to construct limit objects, which will build the basis for our convergence results, we still require some information on the regularity of the time derivatives of the approximate solutions.
\begin{lemma}
\label{lem:time-reg}
Assume $\alpha>1+\sqrt{2}$ and $\beta>1$ and suppose that $\rho:=\min\{\alpha,\beta\}$ satisfies $\rho>\frac{\alpha+1}{\alpha-1}$.
For every $T>0$ exists $C(T)>0$ such that for any $\epsi\in(0,1)$ the solution $(\uep,\vep,\wep)$ of \eqref{approxprob} fulfills
\begin{align}\label{eq:timereg-eq}
\intoT\|u_{\epsi t}\|_{(\W[2,2])^*}+\intoT\big\|\partial_t\ln(\vep+1)\big\|_{(\W[2,2])^*}+\intoT\|w_{\epsi t}\|_{(\W[2,2])^*}\leq C(T).
\end{align}
\end{lemma}

\begin{bew}
For fixed $\psi\in\CSp{\infty}{\bomega}$ we test the first equation of \eqref{approxprob} against $\psi$ and estimate by means of the Hölder and Young inequalities that for all $\epsi\in(0,1)$
\begin{align*}
\Big|\intomega u_{\epsi t}\cdot\psi\Big|&=\Big|-\intomega\nabla\uep\cdot\nabla\psi+\intomega\uep\nabla\wep\cdot\nabla\psi+\intomega f(\uep)\psi\Big|\\
&\leq \Big(\|\nabla\uep\|_{\Lo[2]}^2+\|\uep\nabla\wep\|_{\Lo[2]}^2+\frac{1}{2}\Big)\|\nabla\psi\|_{\Lo[2]}+\|f(\uep)\|_{\Lo[1]}\|\psi\|_{\Lo[\infty]}
\end{align*}
is valid on $(0,\infty)$. Making use of the embedding $\W[2,2]\hookrightarrow\Lo[\infty]$ we can hence find $C_1>0$ such that for all $\epsi\in(0,1)$
\begin{align*}
\int_0^T\!\|u_{\epsi t}\|_{(\W[2,2])^*}\leq C_1\int_0^T\!\Big(\|\nabla\uep\|_{\Lo[2]}^2+\|\uep\nabla\wep\|_{\Lo[2]}^2+\|f(\uep)\|_{\Lo[1]}+\frac{1}{2}\Big),
\end{align*}
which entails the boundedness of the first quantity appearing in \eqref{eq:timereg-eq} in light of Corollary \ref{cor:st-bounds-uep-l2} and Lemmas \ref{lem:vitali-u-nabw} and \ref{lem:equiint-f-g}. Arguing similarly for the second component of \eqref{approxprob} we again fix $\psi\in\CSp{\infty}{\bomega}$ and find that for all $\epsi\in(0,1)$
\begin{align*}
&\Big|\intomega\psi\,\partial_t\ln(\vep+1)\Big|\\=\ &\Big|\intomega\frac{|\nabla\vep|^2}{(\vep+1)^2}\psi-\intomega\frac{\nabla\vep\cdot\nabla\psi}{\vep+1}-\intomega\frac{\vep(\nabla\uep\cdot\nabla\vep)}{(\vep+1)^2}\psi+\intomega\frac{\vep(\nabla\uep\cdot\nabla\psi)}{\vep+1}+\intomega\frac{g(\vep)}{\vep+1}\psi\Big|\\
\leq\ &\Big(\frac{5}{4}\intomega\frac{|\nabla\vep|^2}{(\vep+1)^2}+\intomega|\nabla\uep|^2+\intomega\Big|\frac{g(\vep)}{\vep+1}\Big|\Big)\|\psi\|_{\Lo[\infty]}+\Big(\intomega\frac{|\nabla\vep|^2}{(\vep+1)^2}+\intomega|\nabla\uep|^2+\frac{1}{2}\Big)\|\nabla\psi\|_{\Lo[2]}
\end{align*}
holds on $(0,\infty)$. And that hence, there is $C_2>0$ such that for any $T>0$ and all $\epsi\in(0,1)$
\begin{align*}
\int_0^T\big\|\partial_t\ln(\vep+1)\big\|_{(\W[2,2])^*}\leq C_2\int_0^T\!\Big(\intomega\frac{|\nabla\vep|^2}{(\vep+1)^2}+\intomega|\nabla\uep|^2+\intomega\Big|\frac{g(\vep)}{\vep+1}\Big|+1\Big),
\end{align*}
from which the boundedness of the second part of \eqref{eq:timereg-eq} follows as a consequence of Lemma \ref{lem:st-bounds-vep}, Corollary \ref{cor:st-bounds-uep-l2} and Lemma \ref{lem:equiint-f-g}. Finally, we reiterate the arguments once more for the last equations of \eqref{approxprob}, where for fixed $\psi\in\CSp{\infty}{\bomega}$ and all $\epsi\in(0,1)$, by Hölder's and Young's inequalities and Lemma \ref{lem:Linfty-est-w},
\begin{align*}
\Big|\intomega w_{\epsi t}\psi\Big|&=\Big|-\intomega\nabla\wep\cdot\nabla\psi-\intomega\frac{(\uep+\vep)\wep\psi}{1+\epsi(\uep+\vep)\wep}-\mu\intomega\wep\psi+\intomega r\psi\Big|\\
&\leq \Big(\intomega|\nabla\wep|^2+\frac{1}{4}\Big)\|\nabla\psi\|_{\Lo[2]}+\Big(\ws\intomega\uep+\ws\intomega\vep+\mu\ws|\Omega|+\rs|\Omega|\Big)\|\psi\|_{\Lo[\infty]}
\end{align*}
is valid on $(0,\infty)$. Hence, the embedding $\W[2,2]\hookrightarrow\Lo[\infty]$ readily entails the existence of $C_3>0$ such that for any $T>0$ and all $\epsi\in(0,1)$ we have
\begin{align*}
\int_0^T\big\|w_{\epsi t}\|_{(\W[2,2])^*}\leq C_3\int_0^T\Big(\intomega|\nabla\wep|^2+\intomega\uep+\intomega\vep+1\Big),
\end{align*}
finalizing the proof in light of Lemmas \ref{lem:vitali-nab-w} and \ref{lem:L1-est}.
\end{bew}

\setcounter{equation}{0} 
\section{Construction of limit functions}\label{sec6:limit}
With the $\epsi$-uniform bounds prepared in the previous sections we can now derive the existence of limit functions $u,v,w$ satisfying the regularity conditions appearing in Definition \ref{def:sol}. These uniform bounds, moreover, provide various precompactness properties from which we obtain a list of convergences, which will be of importance when proceeding to prove the inequality and equalities prescribed by \eqref{eq:def-ineq-v}, \eqref{eq:def-eq-u} and \eqref{eq:def-eq-w}.

\begin{proposition}
\label{prop:conv}
Assume $\alpha>1+\sqrt{2}$ and $\beta>1$ and suppose that $\rho:=\min\{\alpha,\beta\}$ satisfies $\rho>\frac{\alpha+1}{\alpha-1}$. Then there exist a sequence $(\epsi_j)_{j\in N}\subset(0,1)$ with $\epsi_j\searrow0$ as $j\to\infty$ and functions
\begin{align*}
u&\in \LSp{\infty}{(0,\infty);\Lo[1]}\cap\LSploc{\alpha+1}{\bomega\times[0,\infty)}\quad\text{with}\quad\nabla u\in\LSploc{2}{\bomega\times[0,\infty);\R^2},\\
v&\in \LSp{\infty}{(0,\infty);\Lo[1]}\cap\LSploc{\beta}{\bomega\times[0,\infty)}\quad\text{with}\quad\nabla\ln(v+1)\in\LSploc{2}{\bomega\times[0,\infty);\R^2},\\
w&\in \LSp{\infty}{\Omega\times(0,\infty)}\cap\LSploc{\rho}{[0,\infty);\W[2,\rho]}\quad\text{with}\quad\nabla w\in\LSploc{2}{\bomega\times[0,\infty);\R^2},
\end{align*}
such that 
\begin{align*}
f(u)\in\LSploc{\frac{\alpha+1}{\alpha}}{\bomega\times[0,\infty)}\quad\text{and}\quad\frac{g(v)}{v+1}\in\LSploc{\frac{\beta}{\beta-1}}{\bomega\times[0,\infty)}
\end{align*}
and such that the solutions $(\uep,\vep,\wep)$ of \eqref{approxprob} fulfill
\begin{alignat}{2}
\uep\to\ & u&&\text{in }\LSploc{p}{\bomega\times[0,\infty)}\text{ for any }p\in[1,\alpha+1)\text{ and a.e. in}\ \Omega\times(0,\infty),\label{eq:conv-u}\\
\uep(\cdot,t)\to\  &u(\cdot,t)&&\text{in }\LSp{2}{\Omega}\text{ for a.e. }t>0,\label{eq:conv-u-steklov}\\
\nabla\uep\wto\  &\nabla u&&\text{in }\LSploc{2}{\bomega\times[0,\infty);\R^2},\label{eq:conv-nab-u}\\
\uep\wto\ &  u&&\text{in }\LSploc{\alpha+1}{\bomega\times[0,\infty)},\label{eq:conv-weak-u}\\
f(\uep)\to\ & f(u)&&\text{in }\LSploc{\frac{\alpha+1}{\alpha}}{\bomega\times[0,\infty)},\label{eq:conv-f}\\
\vep\ \to\ & v&&\text{in }\LSploc{p}{\bomega\times[0,\infty)}\text{ for any }p\in[1,\beta)\text{ and a.e. in}\ \Omega\times(0,\infty),\label{eq:conv-v}\\
\ln(\vep+1)\wto\ & \ln(v+1)\ &&\text{in }\LSploc{2}{[0,\infty);\W[1,2]},\label{eq:conv-nab-v}\\
\vep\wto\ & v&&\text{in }\LSploc{\beta}{\bomega\times[0,\infty)},\label{eq:conv-weak-v}\\
\frac{g(\vep)}{\vep+1}\to\ &\frac{g(v)}{v+1}\ &&\text{in }\LSploc{1}{\bomega\times[0,\infty)},\label{eq:conv-g}\\
\wep\to\ & w&&\text{in }\LSploc{2}{\bomega\times[0,\infty)}\text{ and a.e. in}\ \Omega\times(0,\infty),\label{eq:conv-w}\\
\wep\wsto\ & w&&\text{in }\LSp{\infty}{\Omega\times(0,\infty)},\label{eq:conv-wstar-w}\\
\wep\wto\ & w&&\text{in }\LSploc{\rho}{[0,\infty);\W[2,\rho]},\label{eq:conv-weak-w-W2rho}\\
\nabla \wep \to\ &\nabla w &&\text{in }\LSploc{2}{\bomega\times[0,\infty);\R^2},\label{eq:conv-strong-nab-w}\\
\frac{(\uep+\vep)\wep}{1+\epsi(\uep+\vep)\wep}&\to (u+v)w\ \, &&\text{in }\LSploc{1}{\bomega\times[0,\infty)},\label{eq:conv-uvw}\\
\uep\nabla\wep\to\ & u\nabla w&&\text{in }\LSploc{2}{\bomega\times[0,\infty);\R^2},\label{eq:conv-u-nab-w}
\end{alignat}
as $\epsi=\epsi_j\searrow0$, and such that $u\geq0$, $v\geq0$ and $w\geq0$ a.e. in $\Omega\times(0,\infty)$. Moreover, $u$ and $v$ satisfy
\begin{align}\label{eq:L1-bound-limits}
\intomega u(\cdot,t)\leq \ubs\quad\text{and}\quad\intomega v(\cdot,t)\leq \vbs\quad\text{for a.e. }t>0,
\end{align}
as well as
\begin{align}\label{eq:L1-bound-limits-with-g}
\intomega v(\cdot,t)\leq \intomega v_0+\intotomega g(v)\quad\text{for a.e. }t>0.
\end{align}
\end{proposition}

\begin{bew}
As a consequence of Corollary \ref{cor:st-bounds-uep-l2} and Lemma \ref{lem:time-reg} $\{\uep\}_{\epsi\in(0,1)}$ is bounded in $\LSp{2}{(0,T);\W[1,2]}$ with $\{u_{\epsi t}\}_{\epsi\in(0,1)}$ being bounded in $\LSp{1}{(0,T);(\W[2,2])^*}$ for any $T>0$. Thus, we may employ an Aubin--Lions type lemma (e.g. \cite[Corollary 8.4]{Sim87}) and obtain that $\{\uep\}_{\epsi\in(0,1)}$ is relatively precompact in $\LSp{2}{(0,T);\W[1,2]}$ for any $T>0$. This ensures the existence of $u\in\LSploc{2}{[0,\infty);\W[1,2]}$ and a subsequence $(\epsi_j)_{j\in\N}\subset(0,1)$ satisfying $\uep\to u$ a.e. in $\Omega\times(0,\infty)$ and in $\LSploc{2}{\bomega\times[0,\infty)}$, which also implies \eqref{eq:conv-u-steklov}. The weak convergences (if necessary along another non-relabeled subsequence) featured in \eqref{eq:conv-nab-u} and \eqref{eq:conv-weak-u} are then immediate consequences of the spatio-temporal bounds established in Corollary \ref{cor:st-bounds-uep-l2}. The strong convergences appearing in \eqref{eq:conv-u} and \eqref{eq:conv-f} are then implied by Vitali's theorem, where for the equiintegrability of $\big\{|f(u_{\epsi_j})|^{\frac{\alpha+1}{\alpha}}\big\}_{j\in\N}$ we additionally rely on Lemma \ref{lem:equiint-f-g}. Similarly, a combination of Lemmas \ref{lem:st-bounds-vep} and \ref{lem:time-reg} with a Poincaré-type inequality, the basic estimate $0\leq\ln(\xi+1)\leq \xi$ on $[0,\infty)$ and Lemma \ref{lem:L1-est} entails the boundedness of $\{\ln(\vep+1)\}_{\epsi\in(0,1)}$ in $\LSp{2}{(0,T);\W[1,2]}$ and $\{v_{\epsi t}\}_{\epsi\in(0,1)}$ in $\LSp{1}{(0,T);(\W[2,2])^*}$ for any $T>0$, so that another application of the Aubin--Lions lemma provides a further subsequence along which $\vep\to v$ a.e. in $\Omega\times(0,\infty)$, from which, again due to the boundedness of $\{\ln(\vep+1)\}_{\epsi\in(0,1)}$ in $\LSp{2}{(0,T);\W[1,2]}$, we can also conclude that \eqref{eq:conv-nab-v} holds. The precompactness properties contained in the the bounds provided by Lemmas \ref{lem:spat-temp-alpha-beta} and \ref{lem:equiint-f-g}, the fact that $\frac{\beta}{\beta-1}>1$ and Vitali's theorem then entail the properties listed in \eqref{eq:conv-v}, \eqref{eq:conv-weak-v} and \eqref{eq:conv-g}. Likewise we note that 
Lemmas \ref{lem:Linfty-est-w}, \ref{lem:vitali-nab-w} and \ref{lem:time-reg} imply that $\{\wep\}_{\epsi\in(0,1)}$ is bounded in $\LSp{2}{(0,T);\W[1,2]}$ with $\{w_{\epsi t}\}_{\epsi\in(0,1)}$ being bounded in $\LSp{1}{(0,T);(\W[2,2])^*}$ and that upon another application of the Aubin--Lions lemma we obtain yet another subsequence along which \eqref{eq:conv-w} and $\nabla \wep\wto\nabla w$ in $\LSploc{2}{\bomega\times[0,\infty);\R^2}$ hold. As $\rho>1$ we can than make use of Lemma \ref{lem:vitali-nab-w} and Vitali's theorem to confirm that actually \eqref{eq:conv-strong-nab-w} is valid, whereas the weak convergences \eqref{eq:conv-wstar-w} and \eqref{eq:conv-weak-w-W2rho} result immediately from Lemmas \ref{lem:Linfty-est-w} and \ref{lem:spat-temp-est-W2r-wep}, respectively. Finally, the strong convergence properties in \eqref{eq:conv-uvw} and \eqref{eq:conv-u-nab-w} are a consequence of Vitali's theorem in conjunction with the equiintegrability properties contained in Lemmas \ref{lem:st-bound-l1-uepvepwep} and \ref{lem:vitali-u-nabw}, respectively.

The bounds presented in \eqref{eq:L1-bound-limits} follow from Lemma \ref{lem:L1-est} in light of \eqref{eq:conv-u} and \eqref{eq:conv-v}. To verify \eqref{eq:L1-bound-limits-with-g}, we first note that according to \eqref{eq:conv-v} we can fix a null Set $N\subset(0,\infty)$ such that for all $t\in(0,\infty)\setminus N$ we have $v_\epsi(\cdot,t)\to v(\cdot,t)$ a.e. in $\Omega$ as $\epsi=\epsi_j\searrow0$. Now, relying on Fatou's lemma and the estimate \eqref{eq:spat-temp-minus-plus} from Lemma \ref{lem:spat-temp-alpha-beta}, we find by splitting $g$ into its positive and negative parts that
\begin{align*}
\intomega v(\cdot,t)+\int_0^{t}\!\!\intomega g_-(v)&\leq\liminf_{\epsi=\epsi_j\searrow0}\Big(\intomega \vep(\cdot,t)+\int_0^{t}\!\!\intomega g_-(\vep)\Big)=\liminf_{\epsi=\epsi_j\searrow0}\Big(\intomega v_0+\int_0^{t}\!\!\intomega g_+(\vep)\Big)
\end{align*}
for all $t\in(0,\infty)\setminus N$. Since \eqref{fgprop} implies that $g_+$ is bounded on $(0,\infty)$, we can make use of the dominated convergence theorem and the continuity of $g_+$ to pass to the limit on the right hand side and conclude \eqref{eq:L1-bound-limits-with-g}. The nonnegativity properties are finally an immediate consequence of the nonnegativity of $\uep,\vep$ and $\wep$ established in Lemma \ref{lem:locex} combined with \eqref{eq:conv-u}, \eqref{eq:conv-v} and \eqref{eq:conv-w}, respectively.
\end{bew}

\setcounter{equation}{0} 
{\hypersetup{linkcolor=black}
\section{Solution properties - Proof of Theorem \ref{theo:1}}\label{sec7:solutionprop}}
With the limit objects $(u,v,w)$ provided by the previous proposition satisfying most of the requirements for a generalized solution in the sense of Definition \ref{def:sol}, we are more or less left with checking that \eqref{eq:def-ineq-v} -- \eqref{eq:def-eq-w} are valid for these limit functions. Corresponding identities for the approximating systems \eqref{approxprob} are easily obtainable from straightforward testing methods. It is, however, not clear yet that we may actually pass to the limit in each of the integrals involved. Let us first make sure that the equalities for $u$ and $w$ are satisfied.

\begin{lemma}\label{lem:weak-sol-prop-u-w}
Assume $\alpha>1+\sqrt{2}$ and $\beta>1$ and suppose that $\rho:=\min\{\alpha,\beta\}$ satisfies $\rho>\frac{\alpha+1}{\alpha-1}$. Let $(\epsi_j)_{j\in\N}$ and $(u,v,w)$ be as provided by Proposition \ref{prop:conv}. Then the identities \eqref{eq:def-eq-u} and \eqref{eq:def-eq-w} are satisfied for all $\varphi\in C_0^\infty(\bomega\times[0,\infty))$.
\end{lemma}

\begin{bew}
Testing the first and third equation of \eqref{approxprob} against an arbitrary test function $\varphi\!\in\! C_0^\infty(\bomega\times[0,\infty))$ we obtain
\begin{align*}
-\intinfomega \uep\varphi_t-\intomega u_0\varphi(\cdot,0)=-\intinfomega\nabla\uep\cdot\nabla\varphi+\intinfomega \uep\nabla\wep\cdot\nabla\varphi+\intinfomega f(\uep)\varphi
\end{align*}
and
\begin{align*}
-\!\intinfomega\wep\varphi_t-\!\intomega w_0\varphi(\cdot,0)=-\!\intinfomega\nabla\wep\cdot\nabla\varphi-\!\intinfomega\frac{(\uep+\vep)\wep}{1+\epsi(\uep+\vep)\wep}\varphi-\mu\!\intinfomega\wep\varphi+\!\intinfomega r\varphi,
\end{align*}
respectively, for all $\epsi\in(0,1)$. According to the convergence properties in \eqref{eq:conv-u}, \eqref{eq:conv-nab-u}, \eqref{eq:conv-u-nab-w}, \eqref{eq:conv-f}, \eqref{eq:conv-w}, \eqref{eq:conv-strong-nab-w} and \eqref{eq:conv-uvw} we may pass to the limit in each of the integrals above along the subsequence $(\epsi_j)_{j\in\N}$ provided by Proposition \ref{prop:conv} and readily achieve \eqref{eq:def-eq-u} and \eqref{eq:def-eq-w}.
\end{bew}

The inequality for $v$ appearing in Definition \ref{def:sol} requires more work. In particular, in order to establish an identity corresponding to \eqref{eq:def-ineq-v} we will have to test the second equation of \eqref{approxprob} against $\frac{\psi}{\vep+1}$. The appearing integral $\intoTomega\frac{\vep}{\vep+1}\big(\nabla\ln(\vep+1)\cdot\nabla\uep)\psi$ will then be the greatest adversary for our limit procedure, as it cannot be treated with the convergence results featured in Proposition \ref{prop:conv}. Due to the limited information attainable for $\vep$ a better convergence property for $\nabla\ln(\vep+1)$ appears to be far out of reach, necessitating that we establish a strong convergence for $\nabla\uep$ instead. Thankfully, problems of this kind arise often in chemotaxis systems with generalized solution concepts (e.g. \cite{win15_chemorot,Win-globgensocialint-M3AS18,TB2018_NA}) and can be treated with well-known results for Steklov averages and the convergence properties established in Proprosition \ref{prop:conv}. 

Let us shortly recall that for $T>0$, $p\in[1,\infty]$, $N\geq1$ and $\psi\in\LSploc{p}{\bomega\times(-1,T);\R^N}$ the Steklov averages $S_h\psi\in\LSp{p}{\Omega\times(0,T);\R^N}$, $h\in(0,1)$, are defined by
\begin{align*}
S_h(\psi)(x,t):=\frac{1}{h}\int_{t-h}^t\psi(x,s)\intd s,\quad (x,t)\in\Omega\times(0,T)
\end{align*}
and satisfy $S_h(\psi)\wto\psi$ in $\LSp{p}{\Omega\times(0,T)}$ as $h\searrow0$ if $p\in[1,\infty)$ and $S_h(\psi)\wsto\psi$ in $\LSp{\infty}{\Omega\times(0,T)}$ as $h\searrow0$ if $p=\infty$ (cf. \cite[Lemma 10.2]{win15_chemorot}).

With this notation fixed, we can now enhance arguments akin to those presented in \cite[Lemmas 8.1 and 8.2]{win15_chemorot} and \cite[Lemmas 5.1 and 5.2]{Win-globgensocialint-M3AS18} to cover our current situation, with, in particular, the appearing source term $f(\uep)$ requiring additional attention.

\begin{lemma}
\label{lem:strong-conv-nab-u}
Assume $\alpha>1+\sqrt{2}$ and $\beta>1$ and suppose that $\rho:=\min\{\alpha,\beta\}$ satisfies $\rho>\frac{\alpha+1}{\alpha-1}$. Let $(\epsi_j)_{j\in\N}$ and $u$ be as provided by Proposition \ref{prop:conv}. Then for all $T>0$
\begin{align*}
\nabla\uep\to\nabla u\quad\text{in }\LSp{2}{\Omega\times(0,T);\R^2}\ \text{as}\ \epsi=\epsi_j\searrow 0.
\end{align*}
\end{lemma}

\begin{bew}
Let $T>0$ be given. According to Proposition \ref{prop:conv} we can choose $t_0>T$ such that $t_0$ is a Lebesgue point of the mapping $0<t\mapsto \intomega u^2(\cdot,t)$ and that
\begin{align}\label{eq:steklov-eq1}
\intomega\uep^2(\cdot,t_0)\to\intomega u^2(\cdot,t_0)\quad\text{as }\epsi=\epsi_j\searrow0.
\end{align}
For $\delta\in(0,1)$ and the chosen $t_0>0$ we define the temporal cutoff functions 
\begin{align*}
\zeta_\delta(t)\equiv\zeta^{(t_0)}_\delta(t):=\begin{cases}
1\quad&\text{if }t\in[0,t_0],\\
1-\frac{t-t_0}{\delta}&\text{if }t\in(t_0,t_0+\delta),\\
0&\text{if }t\geq t_0+\delta,
\end{cases}
\end{align*}
which obviously satisfy
\begin{align}\label{eq:steklov-temporal-cutoff-derivative}
\zeta_\delta'\equiv\begin{cases}
-\frac{1}{\delta}\quad&\text{on }(t_0,t_0+\delta),\\
0\quad&\text{on }(0,t_0)\cup(t_0+\delta,\infty).
\end{cases}
\end{align}
Moreover, we fix any sequence $(u_{0k})_{k\in\N}\subset \CSp{1}{\bomega}$ satisfying $u_{0k}\to u_0$ in $\Lo[2]$ as $k\to\infty$ and accordingly let 
\begin{align*}
u_k(x,t):=\begin{cases}
u(x,t)\quad&\text{if }x\in\Omega\text{ and }t>0,\\
u_{0k}(x)\quad&\text{if }x\in\Omega\text{ and }t\leq 0.
\end{cases}
\end{align*}
With these preparations, for $h\in(0,1)$, $\delta\in(0,1)$ and $k\in\N$ we construct the test functions
\begin{align*}
\varphi(x,t):=\zeta_\delta(t) S_h(u_k)(x,t),\quad (x,t)\in\Omega\times(0,\infty),
\end{align*}
which satisfy $\varphi\in\LSp{\infty}{(0,\infty);\W[1,2]}$ with $\varphi_t\in\LSp{2}{\Omega\times(0,\infty)}$ and are compactly supported in $\bomega\times[0,t_0+1]$. Making use of Lemma \ref{lem:weak-sol-prop-u-w} and a completion argument we find that $\varphi$ is an admissible test function for \eqref{eq:def-eq-u}, which entails that
\begin{align}\label{eq:steklov-eq2}
-\intinfomega& u(x,t)\zeta_\delta'(t)S_h (u_k)(x,t)\intd x\intd t\nonumber\\-\int_0^\infty\!&\!\intomega u(x,t)\zeta_\delta(t)\frac{u_k(x,t)-u_k(x,t-h)}{h}\intd x\intd t-\intomega u_0(x) u_{0k}(x)\intd x\\
=&-\intinfomega\nabla u\cdot\zeta_\delta(t)\nabla S_h(u_k)(x,t)\intd x\intd t+\intinfomega u(x,t)\nabla w(x,t)\cdot\zeta_\delta(t)\nabla S_h(u_k)(x,t)\intd x\intd t\nonumber\\
&\quad+\intinfomega f\big(u(x,t)\big)\zeta_\delta(t) S_h(u_k)(x,t)\intd x\intd t\quad\text{for all }h\in(0,1),\;k\in\N\text{ and }\delta\in(0,1).\nonumber
\end{align}
Here, we note that $\nabla u_k\in\LSp{2}{\Omega\times(-1,t_0+1);\R^2}$ and $u_k\in\LSp{\alpha+1}{\Omega\times(-1,t_0+1)}$ entails $\nabla S_h(u_k)=S_h(\nabla u_k)\wto\nabla u_k=\nabla u$ in $\LSp{2}{\Omega\times(0,t_0+1);\R^2}$ and $S_h(u_k)\wto u$ in $\LSp{\alpha+1}{\Omega\times(0,t_0+1)}$ as $h\searrow0$, which also obviously implies that $S_h(u_k)\wto u$ in $\LSp{2}{\Omega\times(0,t_0+1)}$ due to $\alpha>1$. Since Proposition \ref{prop:conv} ensures that $u$ belongs to $\LSp{2}{\Omega\times(0,t_0+1)}$, that $\nabla u$ as well as $u\nabla w$ belong to $\LSp{2}{\Omega\times(0,t_0+1);\R^2}$ and that $f(u)$ belongs to $\LSp{\frac{\alpha+1}{\alpha}}{\Omega\times(0,t_0+1)}$, we can take $h\searrow0$ in each of the integrals on the right hand side and in the first integral on the left. To treat the remaining integral on the left we first employ Young's inequality to obtain that
\begin{align*}
&-\intinfomega u(x,t)\zeta_\delta(t)\frac{u_k(x,t)-u_k(x,t-h)}{h}\\=\ &-\frac{1}{h}\intinfomega\zeta_\delta(t) u_k^2(x,t)\intd x\intd t+\frac{1}{h}\intinfomega \zeta_\delta(t)u_k(x,t) u_k(x,t-h)\intd x\intd t\\
\leq\ &-\frac{1}{2h}\intinfomega \zeta_\delta(t) u_k^2(x,t)\intd x\intd t+\frac{1}{2h}\intinfomega\zeta_\delta(t) u_k^2(x,t-h)\intd x\intd t\\
=\ &-\frac{1}{2h}\intinfomega\zeta_\delta(t) u_k^2(x,t)\intd x\intd t+\frac{1}{2h}\int_{-h}^\infty\!\intomega\zeta_\delta(s+h)u_k^2(x,s)\intd x\intd s\\
=\ &\frac{1}{2}\intinfomega u^2(x,s)\frac{\zeta_\delta(s+h)-\zeta_\delta(s)}{h}\intd x\intd s+\frac{1}{2}\intomega u_{0k}^2(x)\intd x,
\end{align*}
where, by the dominated convergence theorem, we have
\begin{align*}
\frac{1}{2}\intinfomega u^2(x,s)\frac{\zeta_\delta(s+h)-\zeta_\delta(s)}{h}\intd x\intd s\to\frac{1}{2}\intinfomega u^2(x,s)\zeta_\delta'(s)\intd x\intd s\quad\text{as }h\searrow0.
\end{align*}
Hence, passing to the limit $h\searrow0$ in \eqref{eq:steklov-eq2} we can make use of \eqref{eq:steklov-temporal-cutoff-derivative} to obtain that
\begin{align*}
\frac{1}{2\delta}\int_{t_0}^{t_0+\delta}\!&\!\intomega u^2(x,t)\intd x\intd t+\frac{1}{2}\intomega u_{0k}^2(x)\intd x-\intomega u_0(x)u_{0k}(x)\intd x\\
\geq&-\intinfomega\zeta_\delta(t)|\nabla u|^2\intd x\intd t+\intinfomega\zeta_\delta(t) u(x,t)\big(\nabla w(x,t)\cdot\nabla u(x,t)\big)\intd x\intd t\\
&\quad+\intinfomega \zeta_\delta(t)f\big(u(x,t)\big)u(x,t)\intd x\intd t\quad\text{for all }k\in\N\text{ and }\delta\in(0,1).
\end{align*}
Here, as $u_{0k}\to u_0$ in $\Lo[2]$ as $k\to\infty$, we can first let $k\to\infty$ and then make use of the Lebesgue point property of $t_0$ and the dominated convergence theorem to let $\delta\searrow0$ to find, upon appropriate reordering of the terms, that
\begin{align}\label{eq:steklov-eq3}
\int_0^{t_0}\!\!\intomega|\nabla u|^2\geq-\frac{1}{2}\intomega u^2(\cdot,t_0)+\frac{1}{2}\intomega u_0^2+\int_0^{t_0}\!\!\intomega u\nabla u\cdot\nabla w+\int_0^{t_0}\!\!\intomega f(u)u.
\end{align}
In order to link this with an suitable expression concerning the approximate systems \eqref{approxprob} we make use of the strong convergence properties listed in \eqref{eq:conv-f} and \eqref{eq:conv-u-nab-w} of Proposition \ref{prop:conv}, which together with \eqref{eq:conv-weak-u} and \eqref{eq:conv-nab-u}, respectively, yield that
\begin{align*}
&\int_0^{t_0}\!\!\intomega f(\uep)\uep \to\int_0^{t_0}\!\!\intomega f(u) u\quad
\text{ and }\quad\int_0^{t_0}\!\!\intomega \uep\nabla\uep\cdot\nabla\wep \to \int_0^{t_0}\!\!\intomega u\nabla u \cdot\nabla w\quad\text{as }\epsi=\epsi_j\searrow0,
\end{align*}
and conclude from \eqref{eq:steklov-eq3}, and \eqref{eq:steklov-eq1} and from testing the first equation of \eqref{approxprob} against $\uep$ that
\begin{align*}
\int_0^{t_0}\!\!\intomega|\nabla u|^2&\geq\lim_{\epsi=\epsi_j\searrow0}\Big(-\frac{1}{2}\intomega \uep^2(\cdot,t_0)+\frac{1}{2}\intomega u_0^2+\int_0^{t_0}\!\!\intomega \uep\nabla \uep\cdot\nabla \wep+\int_0^{t_0}\!\!\intomega f(\uep)\uep\Big)\\&=\lim_{\epsi=\epsi_j\searrow0}\int_0^{t_0}\!\!\intomega|\nabla\uep|^2.
\end{align*}
Due to $t_0>T$ and Proposition \ref{prop:conv} saying that $\nabla\uep\wto\nabla u$ in $\LSp{2}{\Omega\times(0,t_0);\R^2}$ as $\epsi=\epsi_j\searrow0$, we actually have $\nabla\uep\to\nabla u$ in $\LSp{2}{\Omega\times(0,T);\R^2}$ as $\epsi=\epsi_j\searrow0$, completing the proof.
\end{bew}

With the strong convergence of $\nabla \uep$ in $\LSp{2}{\Omega\times(0,T);\R^2}$ ensured by the previous lemma and the information contained in Proposition \ref{prop:conv}, we have now prepared all the major parts of the proof of Theorem \ref{theo:1} and can verify the solution properties of the limit objects in a straightforward fashion.

\begin{proof}[\textbf{Proof of Theorem \ref{theo:1}}:]
Proposition \ref{prop:conv} provides us with $(u,v,w)$ satisfying all the regularity and nonnegativity conditions required by Definition \ref{def:sol} as well as \eqref{eq:def-mass-ineq} and Lemma \ref{lem:weak-sol-prop-u-w} verifies that \eqref{eq:def-eq-u} and \eqref{eq:def-eq-w} hold, so that we are left with establishing \eqref{eq:def-ineq-v}. We first pick an arbitrary $\psi\in C_0^\infty(\bomega\times[0,\infty))$ and then fix $T>0$ such that $\psi\equiv 0$ in $\Omega\times(T,\infty)$. Testing the second equation of \eqref{approxprob} against $\frac{\psi}{\vep+1}$ yields
\begin{align}\label{eq:theo-proof-eq1}
-\intoTomega\ln(\vep&+1)\psi_t-\intomega\ln(v_0+1)\psi(\cdot,0)\nonumber\\
&=-\intoTomega\nabla\ln(\vep+1)\cdot\nabla\psi+\intoTomega\big|\nabla\ln(\vep+1)\big|^2\psi+\intoTomega\frac{\vep}{\vep+1}(\nabla\uep\cdot\nabla\psi)\\
&\qquad\qquad-\intoTomega\frac{\vep}{\vep+1}\big(\nabla\ln(\vep+1)\cdot\nabla\uep\big)\psi+\intoTomega\frac{g(\vep)}{\vep+1}\psi\nonumber
\end{align}
for all $\epsi\in(0,1)$. Now, with $(\epsi_j)_{j\in\N}$ provided by Proposition \ref{prop:conv} we obtain from \eqref{eq:conv-nab-v} and \eqref{eq:conv-g} that
\begin{align*}
-\intoTomega\ln(\vep+1)\psi_t\to-&\intoTomega\ln(v+1)\psi_t,\quad-\intoTomega\nabla\ln(\vep+1)\cdot\nabla\psi\to-\intoTomega\nabla\ln(v+1)\cdot\nabla\psi\\
\intertext{and}&\intoTomega\frac{g(\vep)}{\vep+1}\psi\to\intoTomega\frac{g(v)}{v+1}\psi\quad\text{as }\epsi=\epsi_j\searrow0.
\end{align*}
On the other hand, we can combine the boundedness of $\big|\frac{\vep}{\vep+1}\big|$ with the fact that $\frac{\vep}{\vep+1}\to\frac{v}{v+1}$ a.e. in $\Omega\times(0,T)$ as $\epsi=\epsi_j\searrow0$ and the strong convergence from Lemma \ref{lem:strong-conv-nab-u} to obtain that also 
$$\frac{\vep}{\vep+1}\nabla\uep\to\frac{v}{v+1}\nabla u\quad\text{in }\LSp{2}{\Omega\times(0,T);\R^2}$$
(see \cite[Lemma 10.4]{win15_chemorot}), so that 
\begin{align*}
\intoTomega\frac{\vep}{\vep+1}(\nabla\uep\cdot\nabla\psi)&\to\intoTomega\frac{v}{v+1}(\nabla u\cdot\nabla\psi)\intertext{and}
-\intoTomega\frac{\vep}{\vep+1}\big(\nabla\ln(\vep+1)\cdot\nabla\uep\big)&\to-\intoTomega\frac{v}{v+1}\big(\nabla\ln(v+1)\cdot\nabla u\big)\quad\text{as }\epsi=\epsi_j\searrow0,
\end{align*}
where the second part is again due to \eqref{eq:conv-nab-v}. Finally, by the lower semicontinuity of the norm in $\LSp{2}{\Omega\times(0,T)}$ with respect to weak convergence we obtain from \eqref{eq:conv-nab-v} that
\begin{align*}
\liminf_{\epsi=\epsi_j\searrow0}\intoTomega\big|\nabla\ln(\vep+1)\big|^2\psi\geq\intoTomega\big|\nabla\ln(v+1)\big|^2\psi.
\end{align*}
Hence, we can pass to the limit in all of the integrals appearing in \eqref{eq:theo-proof-eq1}, which in light of the fact that $\psi\equiv0$ on $\Omega\times(T,\infty)$, entails
\begin{align*}
-\intinfomega\ln(v&+1)\psi_t-\intomega\ln(v_0+1)\psi\\
&\geq-\intinfomega\ln(v+1)\cdot\nabla\psi+\intinfomega\big|\nabla\ln(v+1)\big|^2\psi+\intinfomega\frac{v}{v+1}(\nabla u\cdot\nabla\psi)\\
&\qquad\qquad-\intinfomega\frac{v}{v+1}\big(\nabla\ln(v+1)\cdot\nabla\psi\big)+\intinfomega\frac{g(v)}{v+1}\psi
\end{align*}
and thereby completes the proof.
\end{proof}

\setcounter{equation}{0} 
\section{Conditional eventual regularity estimates}\label{sec8:ev-cond-sob-reg}

In order to prove the eventual smoothing process postulated by Theorem \ref{theo:2} we will have to significantly improve our knowledge on the regularity of $\uep$ and, especially, $\vep$. In this section we will hence take an in-depth look at eventual regularity estimates which can be derived under the additional assumptions that at least beyond some waiting time $T>0$ we have certain (spatio-temporal) bounds at hand for $\uep$, $\vep$ and $\wep$. We begin with a slightly refined version of Lemma \ref{lem:spat-temp-est-W2r-wep}.

\begin{corollary}
\label{cor:cond-spat-temp-est-W2q-wep}
Suppose that for some $q>1$ there are $T\geq0$ and $K>0$ such that for any $\epsi\in(0,1)$ and all $t>T$ the solution $(\uep,\vep,\wep)$ of \eqref{approxprob} satisfies 
\begin{align}\label{eq:cond-spat-temp-est-W2q-wep}
\int_t^{t+1}\!\!\intomega\uep^q+\int_t^{t+1}\!\!\intomega\vep^q \leq K.
\end{align}
Then, there exists $C>0$ such that for any $\epsi\in(0,1)$ and all $t> T+1$
\begin{align*}
\int_t^{t+1}\!\|\wep\|_{\W[2,q]}^q\leq C.
\end{align*}
\end{corollary}

\begin{bew}
This result is quite close to the one presented in Lemma \ref{lem:spat-temp-est-W2r-wep}. The main difference is the possibility of \eqref{eq:cond-spat-temp-est-W2q-wep} only being valid after a waiting time $T\geq0$ for which we pay by the fact that the bound on $\wep$ only holds true for $t>T+1$. The proof, however, works along similar lines as Lemma \ref{lem:spat-temp-est-W2r-wep}. We only recall that upon picking $\tau>T$, letting $\xi:=\xi_\tau$ and employing the maximal Sobolev regularity estimates to the PDE corresponding to $\xi\wep$ we find that there is some $C_1>0$ such that
\begin{align*}\int_\tau^{\tau+2}\!\|\xi\wep\|_{\W[2,q]}^q
&\leq C_1+C_1\int_\tau^{\tau+2}\!\!\intomega\uep^q+C_1\int_\tau^{\tau+2}\!\!\intomega\vep^q\leq C_1+4K C_1,
\end{align*}
which, due to the arbitrary choice of $\tau>T$ and the fact that $\xi\equiv 1$ on $(\tau+1,\tau+2)$, implies that there is some $C_2>0$ such that
\begin{align*}
\int_t^{t+1}\!\|\wep\|_{\W[2,q]}^q\leq C_2
\end{align*} 
holds for each $\epsi\in(0,1)$ and all $t>T+1$.
\end{bew}

If we prescribe even stronger known eventual boundedness properties we can significantly improve the regularity of $\uep$ and $\vep$ in the next result. Due to its cyclic nature this result even provides a basis for an iterative argument we will later employ in Section \ref{sec11:ev-smooth}.

\begin{lemma}
\label{lem:ev-max-sob}
Assume $\alpha>1+\sqrt{2}$ and $\beta>2$. Let $q>\max\{\beta-1,\frac{1}{\beta-2}\}$. If there are $\tau_0>0$ and $K_0>0$ such that
\begin{align}\label{eq:lem-ev-max-sob-wep-vep-bound}
\int_t^{t+1}\!\big\|\wep\big\|_{\W[2,q]}^q+\int_t^{t+1}\!\!\intomega\vep^{q-\beta+2}\leq K_0\quad\text{for all }t>\tau_0\text{ and }\epsi\in(0,1),
\end{align}
and if for all $p\in[1,\infty)$ there are $\tau_p>0$ and $K_p>0$ such that
\begin{align}\label{eq:lem-ev-max-sob-uep-bound}
\intomega\uep^p(\cdot,t)\leq K_p\quad\text{for all }t>{\tau}_p\text{ and }\epsi\in(0,1),
\end{align}
then for all $\theta\in(1,q)$ there exist $\Ts>0$ and $C>0$ such that for each $\epsi\in(0,1)$ and all $t>\Ts$
\begin{align}\label{eq:lem-ev-sob-uep}
\int_t^{t+1}\!\big\|\uep\big\|_{\W[2,\theta]}^\theta+\int_t^{t+1}\!\!\intomega|\nabla\uep|^{2\theta}\leq C
\end{align}
and, moreover, there exists $C'>0$ such that for each $\epsi\in(0,1)$ and all $t>\Ts$
\begin{align*}
\intomega\vep^{q-\beta+2}(\cdot,t)+\int_t^{t+1}\!\!\intomega\vep^{q+1}+\int_t^{t+1}\!\!\intomega\vep^{q-\beta}|\nabla\vep|^2\leq C'\qquad\text{and}\qquad\int_t^{t+1}\!\big\|\wep\big\|_{\W[2,q+1]}^{q+1}\leq C'.
\end{align*}
\end{lemma}

\begin{bew}
We pick any $\theta\in(1,q)$ and let $a:=\max\big\{\alpha\theta,\frac{q\theta}{q-\theta}\big\}>1$. According to \eqref{eq:lem-ev-max-sob-wep-vep-bound} and \eqref{eq:lem-ev-max-sob-uep-bound} for $\tau':=\max\{\tau_0,\tau_a\}>0$ and $K':=2\max\{K_0, K_a\}>0$ we have
\begin{align}\label{eq:ev-max-sob-eq1}
\int_t^{t+1}\!\intomega\big|\Delta\wep\big|^q+\int_t^{t+1}\!\intomega\uep^a\leq K'\quad\text{for any }\epsi\in(0,1)\text{ and all }t>\tau'.
\end{align} 
For $t_0>\tau'$ we denote by $\xi:=\xi_{t_0}$ the cut-off function provided by Definition \ref{def:cutoff} and observe that $\xi\uep$ solves the equation
\begin{align}\label{eq:ev-max-sob-eqcut}
\begin{split}
\big(\xi\uep\big)_t&=\Delta(\xi\uep)-\nabla(\xi\uep)\cdot\nabla\wep-\xi\uep\Delta\wep+\xi f(\uep)+\xi'\uep\quad\text{in }\Omega\times(t_0,\infty)\\
&\text{with }\big(\xi\uep\big)(\cdot,t_0)=0\quad\text{in }\Omega\quad\text{and}\quad \nabla(\xi\uep)\cdot\nu=0\quad\text{on }\romega\times(t_0,\infty).
\end{split}
\end{align}
Now we once more rely on the maximal Sobolev regularity estimates (\cite{HieberPruss-CPDE92,GigSohr91}) and apply them to \eqref{eq:ev-max-sob-eqcut}. Thus, we find that for any $\theta\in(1,q)$ there is $C_1>0$ such that
\begin{align*}
\int_{t_0}^{t_0+2}\!\big\| (\xi\uep)_t\big\|_{\Lo[\theta]}^\theta+\int_{t_0}^{t_0+2}\!\big\|\xi\uep\big\|_{\W[2,\theta]}^\theta&\leq C_1\int_{t_0}^{t_0+2}\!\big\|\nabla(\xi\uep)\cdot\nabla\wep\big\|_{\Lo[\theta]}^\theta+C_1\int_{t_0}^{t_0+2}\!\big\|\xi\uep\Delta\wep\big\|_{\Lo[\theta]}^\theta\\&\quad+C_1\int_{t_0}^{t_0+2}\!\big\|\xi f(\uep)\big\|_{\Lo[\theta]}^\theta+C_1\int_{t_0}^{t_0+2}\big\|\xi'\uep\big\|_{\Lo[\theta]}^\theta
\end{align*}
for all $\epsi\in(0,1)$. Making use of the fact that $\theta<q$, Young's inequality and \eqref{fgprop} yield $C_2>0$ such that
\begin{align*}
&\int_{t_0}^{t_0+2}\!\big\|\xi\uep\big\|_{\W[2,\theta]}^\theta\\
\leq\ &C_1\int_{t_0}^{t_0+2}\!\!\intomega\big|\nabla(\xi\uep)\big|^{\frac{2q\theta}{2q-\theta}}+C_1\int_{t_0}^{t_0+2}\!\!\intomega\big|\nabla\wep\big|^{2q}+C_1\int_{t_0}^{t_0+2}\!\!\intomega\big|\Delta\wep\big|^q+C_1\int_{t_0}^{t_0+2}\!\!\intomega\uep^\frac{q\theta}{q-\theta}\\&\quad+C_2\int_{t_0}^{t_0+2}\!\!\intomega\uep^{\alpha\theta}+C_2+C_1\|\xi'\|_{\LSp{\infty}{\R}}^\theta\int_{t_0}^{t_0+2}\!\intomega\uep^\theta\quad\text{for all }\epsi\in(0,1),
\end{align*}
where we also used that $\xi\leq 1$ on $\R$. To further estimate the terms on the right hand side, we note that in light of the \GNI\ (e.g. \cite[Theorem I.9.3]{fr69} combined with \cite[Theorem I.19.1]{fr69}) for any $1<s\leq\sigma$ we can find $C_{s,\sigma}>0$ such that
\begin{align}\label{eq:ev-max-sob-gni}
\|\nabla\phi\|_{\Lo[2s]}^{2s}\leq C_{s,\sigma}\|\phi\|_{\W[2,\sigma]}^s\|\phi\|_{\Lo[\frac{s\sigma}{\sigma-s}]}^s\quad\text{for all }\phi\in\W[2,\sigma].
\end{align}
Since $\theta<q$ implies $\frac{q\theta}{2q-\theta}<\theta$, we can employ \eqref{eq:ev-max-sob-gni} once for $\sigma=\theta$ and $s=\frac{q\theta}{2q-\theta}$ and once for $\sigma=s=q$ to conclude the existence of $C_3
>0$ satisfying
\begin{align*}
&\int_{t_0}^{t_0+2}\!\big\|\xi\uep\big\|_{\W[2,\theta]}^\theta\nonumber\\
\leq\ &C_1C_3\int_{t_0}^{t_0+2}\!\big\|\xi\uep\big\|_{\W[2,\theta]}^{\frac{q\theta}{2q-\theta}}\big\|\xi\uep\big\|_{\Lo[\frac{q\theta}{q-\theta}]}^\frac{q\theta}{2q-\theta}+C_1C_3\int_{t_0}^{t_0+2}\!\big\|\wep\big\|_{\W[2,q]}^q\big\|\wep\big\|_{\Lo[\infty]}^q+C_1\int_{t_0}^{t_0+2}\!\!\intomega\big|\Delta\wep\big|^q\nonumber\\
&\quad+C_1\int_{t_0}^{t_0+2}\!\!\intomega\uep^\frac{q\theta}{q-\theta}+C_2\int_{t_0}^{t_0+2}\!\!\intomega\uep^{\alpha\theta}+C_2+C_1\|\xi'\|_{\LSp{\infty}{\R}}^\theta\int_{t_0}^{t_0+2}\!\intomega\uep^\theta\quad\text{for all }\epsi\in(0,1).
\end{align*}
Drawing once more on the fact that $\theta<q$ implies $\frac{q\theta}{2q-\theta}<\theta$, we can rely on Young's inequality,  Lemma \ref{lem:Linfty-est-w},  and $\xi\leq 1$ on $\R$ to find $C_4>0$ such that
\begin{align*}
&\int_{t_0}^{t_0+2}\!\big\|\xi\uep\big\|_{\W[2,\theta]}^\theta\nonumber\\\leq\ &\frac{1}{2}\int_{t_0}^{t_0+2}\!\big\|\xi\uep\big\|_{\W[2,\theta]}^\theta+C_4\int_{t_0}^{t_0+2}\!\intomega\uep^\frac{q\theta}{q-\theta}+C_1C_3\ws^{q}\int_{t_0}^{t_0+2}\!\big\|\wep\big\|_{\W[2,q]}^q+C_1\int_{t_0}^{t_0+2}\!\intomega\big|\Delta\wep\big|^q\\
&\qquad+C_1\int_{t_0}^{t_0+2}\!\intomega\uep^\frac{q\theta}{q-\theta}+C_2\int_{t_0}^{t_0+2}\!\intomega\uep^{\alpha\theta}+C_2+C_1\|\xi'\|_{\Lo[\infty]}^\theta\int_{t_0}^{t_0+2}\!\intomega\uep^\theta\quad\text{for all }\epsi\in(0,1).
\end{align*}
Due to $\alpha>1$ and $a=\max\{\alpha\theta,\frac{q\theta}{q-\theta}\}$, we can apply Young's inequality again, this time in the four integrals only containing $\uep$, to obtain $C_5>0$ such that
\begin{align*}
\frac{1}{2}\int_{t_0}^{t_0+2}\!\big\|\xi\uep\big\|_{\W[2,\theta]}^\theta
\leq\ &C_5\int_{t_0}^{t_0+2}\!\!\intomega\uep^a+C_1C_3\ws^{q}\int_{t_0}^{t_0+2}\big\|\wep\big\|_{\W[2,q]}^q+C_1\int_{t_0}^{t_0+2}\!\intomega\big|\Delta\wep\big|^q+C_5
\end{align*}
is valid for all $\epsi\in(0,1)$. In light of \eqref{eq:ev-max-sob-eq1}, Lemma \ref{lem:Linfty-est-w} combined with \cite[Theorem I.19.1]{fr69}, the arbitrary choice of $t_0>\tau'$ and the fact that $\xi\equiv 1$ on $(t_0+1,t_0+2)$ this entails the existence of $C_6>0$ such that
\begin{align*}
\int_{t}^{t+1}\!\big\|\uep\big\|_{\W[2,\theta]}^\theta\leq C_6\quad\text{for each }\epsi\in(0,1)\text{ and all }t>T_1:=\tau'+1,
\end{align*}
proving the first bound featured in \eqref{eq:lem-ev-sob-uep}. Next, for any $\theta\in(1,q)$ we pick $\theta^*\in(\theta,q)$ and find that, again in light of \eqref{eq:ev-max-sob-gni}, there are $T_2> T_1$ and $C_7>0$ such that
\begin{align*}
\int_t^{t+1}\!\!\intomega\big|\nabla\uep\big|^{2\theta}\leq C_{\theta,\theta^*}\int_t^{t+1}\!\big\|\uep\big\|_{\W[2,\theta^*]}^{\theta}\big\|\uep\big\|_{\Lo[\frac{\theta\theta^*}{\theta^*-\theta}]}^\theta\leq C_7\quad\text{for each }\epsi\in(0,1)\text{ and all }t>T_2,
\end{align*}
due to \eqref{eq:lem-ev-max-sob-uep-bound}, the choice $\theta^*<q$ and the first part of this lemma applied to $\theta^*$ instead of $\theta$. To verify the asserted boundedness for $\vep$, we note that $q>\max\{\beta-1,\frac{1}{\beta-2}\}$ implies that $\lambda:=q-\beta+2>1$. Then, we turn to testing the second equation of \eqref{approxprob} against $\vep^{\lambda-1}$ and find upon integrating by parts that
\begin{align*}
\frac{1}{\lambda}\frac{\intd}{\intd t}\intomega\vep^\lambda(\cdot,t)=-(\lambda-1)\intomega\vep^{\lambda-2}\big|\nabla\vep\big|^2+(\lambda-1)\intomega\vep^{\lambda-1}(\nabla\vep\cdot\nabla\uep)+\intomega\vep^{\lambda-1}g(\vep)
\end{align*}
holds on $(0,\infty)$. Therefore, drawing on the nonnegativity of $\vep$, \eqref{fgprop} and Young's inequality we obtain $C_8>0$ and $C_9>0$ such that for each $\epsi\in(0,1)$
\begin{align}\label{eq:ev-max-sob-eq3}
\frac{1}{\lambda}\frac{\intd}{\intd t}\intomega\vep^\lambda(\cdot,t)&\leq-(\lambda-1)\intomega\vep^{\lambda-2}\big|\nabla\vep\big|^2+(\lambda-1)\intomega\vep^{\lambda-1}(\nabla\vep\cdot\nabla\uep)-K_g\intomega\vep^{\lambda+\beta-1}+L_g\intomega\vep^{\lambda-1}\nonumber\\
&\leq-\frac{\lambda-1}{2}\intomega\vep^{\lambda-2}\big|\nabla\vep\big|^2-\frac{K_g}{2}\intomega\vep^{\lambda+\beta-1}+ \frac{\lambda-1}{2}\intomega\vep^{\lambda}\big|\nabla\uep\big|^{2}+C_8\nonumber\\
&\leq-\frac{\lambda-1}{2}\intomega\vep^{\lambda-2}\big|\nabla\vep\big|^2-\frac{K_g}{4}\intomega\vep^{\lambda+\beta-1}+C_9\intomega\big|\nabla\uep\big|^{\frac{2(\lambda+\beta-1)}{\beta-1}}+C_8
\end{align}
is valid on $(0,\infty)$. Since \eqref{eq:lem-ev-max-sob-wep-vep-bound} implies that for all $t>\tau_0+1$ and any $\epsi\in(0,1)$ there is some $T_t(\epsi)\in(t-1,t)$ such that
\begin{align*}
\intomega\vep^\lambda(\cdot,T_t(\epsi))\leq K_0\quad\text{for any }\epsi\in(0,1),
\end{align*}
we find from \eqref{eq:ev-max-sob-eq3} that for any arbitrary $t>\tau_0+1$ and $\epsi\in(0,1)$
\begin{align*}
\frac{1}{\lambda}\intomega\vep^\lambda(\cdot,t)\leq \frac{1}{\lambda}\intomega\vep^\lambda(\cdot,T_t(\epsi))+C_9\int_{T_t(\epsi)}^t\!\intomega\big|\nabla\uep\big|^{\frac{2(\lambda+\beta-1)}{\beta-1}}+C_8(t-T_t(\epsi)),
\end{align*}
which due to $\frac{\lambda+\beta-1}{\beta-1}=\frac{q+1}{\beta-1}<q$, \eqref{eq:lem-ev-max-sob-uep-bound}, the fact that $t-T_t(\epsi)<1$ for all $\epsi\in(0,1)$ and the arbitrary choice of $t>\tau_0+1$ yields $T_3>\max\{T_2,\tau_0+1\}$ and $C_{10}>0$, both independent of $\epsi\in(0,1)$, such that
\begin{align*}
\frac{1}{\lambda}\intomega\vep^{q-\beta+2}(\cdot,t)\leq \frac{K_0}{\lambda}+C_{10}\quad\text{for each }\epsi\in(0,1)\text{ and all }t>T_3.
\end{align*}
With this, we return to \eqref{eq:ev-max-sob-eq3} and find $C_{11}>0$ such that
\begin{align*}
\frac{K_g}{4}&\int_t^{t+1}\!\!\intomega\vep^{q+1}+\frac{q-\beta+1}{2}\int_t^{t+1}\!\!\intomega\vep^{q-\beta}|\nabla\vep|^2\\\leq\frac{1}{q-\beta+2}&\intomega\vep^{q-\beta+2}(\cdot,t)+C_9\int_t^{t+1}\!\!\intomega\big|\nabla\uep\big|^\frac{2(q+1)}{\beta-1}+C_8\leq C_{11}\quad\text{for each }\epsi\in(0,1)\text{ and all }t>T_3,
\end{align*}
as claimed. Finally, we observe that the bound on $\int_t^{t+1}\!\!\intomega\vep^{q+1}$ for all $t>T_3$, when combined with the fact that by \eqref{eq:lem-ev-max-sob-uep-bound} there is some $T_4\geq T_3$ such that $\uep$ is bounded in $\Lo[q+1]$ for all $t>T_4$ and an application of Corollary \ref{cor:cond-spat-temp-est-W2q-wep}, immediately yields $C_{11}>0$ such that
\begin{align*}
\int_t^{t+1}\!\!\big\|\wep\big\|_{\W[2,q+1]}^{q+1}\leq C_{11}\quad\text{for each }\epsi\in(0,1) \text{ and all }t>T_4,
\end{align*}
which completes the proof.
\end{bew}

\setcounter{equation}{0} 
\section{Decay properties of \texorpdfstring{$\wep$}{w}}\label{sec9:w-decay}
While at least some of the bounds appearing in the conditions of Lemma \ref{lem:ev-max-sob} are similar to those already established in previous sections, we are still missing the required eventual $L^q$ regularity of $\uep$ for arbitrary large $q\geq1$. Comparing with similar chemotaxis systems where the eventual smoothness of solutions could be observed, one approach to this boundedness consists of a combination of a decay property of the nutrient density and a differential inequality for the quantity $\intomega\frac{\uep^q}{(2\delta-\wep)^\theta}$ with small $\delta>0$ and suitably chosen $\theta,T>0$ (see e.g. \cite{Lan-Longterm_M3AS16,win-stab2d-ArchRatMechAna12}). To make use of this approach, however, we first have to ensure that $\wep$ indeed does decay below any given threshold, provided the waiting time is large enough. As an initial step we will ensure that at least $\|\wep\|_{\LSp{1}{\Omega\times(t,t+1)}}$ does decay past any positive threshold.

\begin{lemma}
\label{lem:dec-uep-vep-wep}
Assume $\mu>0$, $\alpha>1+\sqrt{2}$ and $\beta>1$ and suppose that $\rho:=\min\{\alpha,\beta\}$ satisfies $\rho>\frac{\alpha+1}{\alpha-1}$ and that $r\in\CSp{1}{\bomega\times[0,\infty)}\cap\LSp{\infty}{\Omega\times(0,\infty)}\cap\LSp{1}{(0,\infty);\Lo[\infty]}$ fulfills \eqref{rprop2}. Then for any $\delta>0$ there exists $T>0$ such that for all $\epsi\in(0,1)$ the solution $(\uep,\vep,\wep)$ of \eqref{approxprob} satisfies
\begin{align*}
\int_T^{\infty}\!\!\intomega\frac{(\uep+\vep)\wep}{1+\epsi(\uep+\vep)\wep}+\int_T^{\infty}\!\!\intomega \wep< \delta.
\end{align*}
\end{lemma}

\begin{bew}
We infer from the integrating third equation of \eqref{approxprob} and the divergence theorem that
\begin{align*}
\frac{\intd}{\intd t}\intomega\wep(\cdot,t)+\intomega\Big(\frac{(\uep+\vep)\wep}{1+\epsi(\uep+\vep)\wep}\Big)(\cdot,t)+\mu\intomega\wep(\cdot,t)\leq \intomega r(\cdot,t)
\end{align*}
for all $t>0$. Hence, due to the nonnegativity of $\wep$ on $\Omega\times(0,\infty)$ and \eqref{rprop2},
\begin{align*}
\int_0^t\!\!\intomega\frac{(\uep+\vep)\wep}{1+\epsi(\uep+\vep)\wep}+\mu\int_0^t\!\!\intomega\wep\leq\intomega w_0+ \rhs
\end{align*}
is valid for all $t>0$, from which we can readily conclude the desired outcome.
\end{bew}

Following the ideas of \cite[Lemma 6.7]{Win-globgensocialint-M3AS18} we can now verify that also $\|\wep(\cdot,t)\|_{\Lo[\infty]}$ decays to zero.

\begin{lemma}
\label{lem:dec-wep-Linfty}
Assume $\mu>0$, $\alpha>1+\sqrt{2}$ and $\beta>1$ and suppose that $\rho:=\min\{\alpha,\beta\}$ satisfies $\rho>\frac{\alpha+1}{\alpha-1}$ and that $r\in\CSp{1}{\bomega\times[0,\infty)}\cap\LSp{\infty}{\Omega\times(0,\infty)}\cap\LSp{1}{(0,\infty);\Lo[\infty]}$ fulfills \eqref{rprop2}. Then for any $\delta>0$ there is $T>0$ such that for any $t>T$ and all $\epsi\in(0,1)$ the solution $(\uep,\vep,\wep)$ of \eqref{approxprob} satisfies
\begin{align*}
\big\|\wep(\cdot,t)\big\|_{\Lo[\infty]}<\delta.
\end{align*}
\end{lemma}

\begin{bew}
According to well-known estimates for the Neumann heat semigroup (e.g. \cite[Lemma 1.3]{win10jde}) we can find $C_1>0$ such that
\begin{align*}
\|e^{\sigma\Delta}\varphi\|_{\Lo[\infty]}\leq C_1 \sigma^{-1}\|\varphi\|_{\Lo[1]}\quad\text{ for all }\sigma\in(0,1)\text{ and }\varphi\in\CSp{0}{\bomega}\text{ with }\overline{\varphi}=0,
\end{align*}
where $\overline{\varphi}:=\tfrac{1}{|\Omega|}\intomega\varphi$. Now, given any arbitrary $\delta>0$ we note that by Lemma \ref{lem:dec-uep-vep-wep} there is some $T_1>1$ such that for all $t>T_1$ and any $\epsi\in(0,1)$
\begin{align}\label{eq:dec-wep-Linfty-eq1}
2C_1\int_{t-1}^{t-\frac12}\!\|\wep(\cdot,s)\|_{\Lo[1]}\intd s<\frac{\delta}{4},
\end{align}
and that in light of \eqref{rprop2} we can pick $T_2>1$ such that
\begin{align*}
\int_{t-1}^{t}\!\|r(\cdot,s)\|_{\Lo[\infty]}\intd s<\frac{\delta}{2}\quad\text{for all }t>T_2.
\end{align*}
Letting $T_0:=\max\{T_1,T_2\}$  we fix an arbitrary $t>T_0$ and note that from \eqref{eq:dec-wep-Linfty-eq1} we can conclude that there is some $t_0\in(t-1,t-\frac12)$ satisfying
\begin{align*}
2C_1\big\|\wep(\cdot,t_0)\big\|_{\Lo[1]}<\frac{\delta}{2}\quad\text{for all }\epsi\in(0,1).
\end{align*}
With this we now observe that the comparison principle implies that
\begin{align*}
\wep(\cdot,t)&=e^{(t-t_0)\Delta}\wep(\cdot,t_0)-\int_{t_0}^t e^{(t-s)\Delta}\Big(\frac{(\uep+\vep)\wep}{1+\epsi(\uep+\vep)\wep}(\cdot,s)+\mu\wep(\cdot,s)\Big)\intd s+\int_{t_0}^t e^{(t-s)\Delta} r(\cdot,s)\intd s\\
&\leq e^{(t-t_0)\Delta}\big(\wep(\cdot,t_0)-\overline{\wep(\cdot,t_0)}\big)+\overline{\wep(\cdot,t_0)}+\int_{t_0}^t\big\|r(\cdot,s)\big\|_{\Lo[\infty]}\intd s\quad\text{in }\Omega,
\end{align*}
and that hence
\begin{align*}
\big\|\wep(\cdot,t)\big\|_{\Lo[\infty]}&\leq C_1(t-t_0)^{-1}\big\|\wep(\cdot,t_0)\|_{\Lo[1]}+\big\|\overline{\wep(\cdot,t_0)}\big\|_{\Lo[\infty]}+\int_{t_0}^{t}\big\|r(\cdot,s)\big\|_{\Lo[\infty]}\intd s\\
&\leq \big(2C_1+\tfrac{1}{|\Omega|}\big)\big\|\wep(\cdot,t_0)\big\|_{\Lo[1]}+\int_{t_0}^{t_0+1}\big\|r(\cdot,s)\big\|_{\Lo[\infty]}\intd s< \frac{\delta}{2}+\frac{\delta}{2}=\delta
\end{align*}
for all $\epsi\in(0,1)$, completing the proof.
\end{bew}

\setcounter{equation}{0} 
\section{Eventual \texorpdfstring{$L^p$}{Lp} regularity of \texorpdfstring{$\uep$}{u}}\label{sec10:ev-lp-uep}

With the decay property of $\wep$ at hand, we can now turn to tracking the time-evolution of $\intomega\frac{\uep^q}{(2\delta-\wep)^\theta}$. Functionals of this type have successfully been utilized in e.g. \cite[Lemma 3.5]{Lan-Longterm_M3AS16} and \cite[Lemma 5.1]{win-stab2d-ArchRatMechAna12} in similar chemotaxis system.
\begin{lemma}
\label{lem:testing-uep-wep}
Let $q>1$, $\theta>0$, $\delta\in(0,\frac1{2q})$ and $T>0$. Then there is $C>0$ such that whenever the solution $(\uep,\vep,\wep)$ of \eqref{approxprob} fulfills
\begin{align*}
\|\wep(\cdot,t)\|_{\Lo[\infty]}\leq\delta\quad\text{for all }\epsi\in(0,1)\text{ and }t>T,
\end{align*}
then the inequality
\begin{align*}
\frac{\intd}{\intd t}\intomega\frac{\uep^q}{(2\delta-\wep)^\theta}&+qK_f\intomega\frac{\uep^{q+\alpha-1}}{(2\delta-\wep)^\theta}\\
&\leq \bigg[\frac{\big(2q\theta+2q(q-1)\delta\big)^2}{4\big(\theta(\theta+1)-2q\theta\delta\big)}-q(q-1)\bigg]\intomega\frac{\uep^{q-2}|\nabla\uep|^2}{(2\delta-\wep)^\theta}+C\intomega\frac{\uep^q}{(2\delta-\wep)^\theta}+C
\end{align*}
holds for all $t>T$ and $\epsi\in(0,1)$.
\end{lemma}

\begin{bew}
Straightforward calculations utilizing integration by parts shows 
\begin{align*}
&\frac{\intd}{\intd t}\intomega\frac{\uep^q}{(2\delta-\wep)^\theta}\\
=\ &q\!\intomega\frac{\uep^{q-1}\big(\Delta\uep-\nabla\cdot(\uep\nabla\wep)+f(\uep)\big)}{(2\delta-\wep)^\theta}+\theta\!\intomega\frac{\uep^q}{(2\delta-\wep)^{\theta+1}}\big(\Delta\wep-\frac{(\uep+\vep)\wep}{1+\epsi(\uep+\vep)\wep}-\mu\wep+r\big)\\
=\ &-q(q-1)\!\intomega\frac{\uep^{q-2}|\nabla\uep|^2}{(2\delta-\wep)^{\theta}}-2q\theta\!\intomega\frac{\uep^{q-1}}{(2\delta-\wep)^{\theta+1}}(\nabla\uep\cdot\nabla\wep)+q(q-1)\!\intomega\frac{\uep^{q-1}}{(2\delta-\wep)^\theta}(\nabla\uep\cdot\nabla\wep)\\
&\quad +q\theta\!\intomega\frac{\uep^q|\nabla\wep|^2}{(2\delta-\wep)^{\theta+1}}+q\!\intomega\frac{\uep^{q-1} f(\uep)}{(2\delta-\wep)^\theta}-\theta(\theta+1)\!\intomega\frac{\uep^q|\nabla\wep|^2}{(2\delta-\wep)^{\theta+2}}\\
&\qquad-\theta\!\intomega\frac{(\uep+\vep)\wep\uep^q}{(1+\epsi(\uep+\vep)\wep)(2\delta-\wep)^{\theta+1}}-\mu\theta\!\intomega\frac{\uep^q\wep}{(2\delta-\wep)^{\theta+1}}+\theta\!\intomega\frac{r\uep^q}{(2\delta-\wep)^{\theta+1}}\\
\end{align*}
for all $t>T$. Therefore, relying on \eqref{fgprop}, \eqref{rprop}, the nonegativity of $\uep,\vep,\wep$ and $\mu$, as well as on the fact that $\frac{1}{2\delta}\leq\frac{1}{2\delta-\wep}\leq\frac{1}{\delta}$ we can estimate
\begin{align}\label{eq:testing-uep-eq1}
&\frac{\intd}{\intd t}\intomega\frac{\uep^q}{(2\delta-\wep)^\theta}\nonumber\\
\leq\ &-q(q-1)\!\!\intomega\frac{\uep^{q-2}|\nabla\uep|^2}{(2\delta-\wep)^\theta}+\!\big(2q\theta+2q(q-1)\delta\big)\!\!\intomega \frac{\uep^{q-1}|\nabla\uep||\nabla\wep|}{(2\delta-\wep)^{\theta+1}}-\!\big(\theta(\theta+1)-2q\theta\delta\big)\!\!\intomega \frac{\uep^q|\nabla\wep|^2}{(2\delta-\wep)^{\theta+2}}\nonumber\\
&\quad-q K_f\!\!\intomega\frac{\uep^{q+\alpha-1}}{(2\delta-\wep)^\theta}+q L_f\!\!\intomega\frac{\uep^{q-1}}{(2\delta-\wep)^\theta}+\frac{\theta\rs}{\delta}\intomega\frac{\uep^q}{(2\delta-\wep)^{\theta}}
\end{align}
for all $t>T$. To treat some of the terms further, we apply Young's inequality to the second term on the right, yielding 
\begin{align*}
\intomega\frac{\uep^{q-1}|\nabla\uep||\nabla\wep|}{(2\delta-\wep)^{\theta+1}}\leq \frac{\big(2q\theta+2q(q-1)\delta\big)}{4\big(\theta(\theta+1)-2q\theta\delta\big)}\intomega\frac{\uep^{q-2}|\nabla\uep|^2}{(2\delta-\wep)^\theta}+\frac{\big(\theta(\theta+1)-2q\theta\delta\big)}{\big(2q\theta+2q(q-1)\delta\big)}\intomega\frac{\uep^q|\nabla\wep|^2}{(2\delta-\wep)^{\theta+2}}
\end{align*}
and to the second to last term on the right, providing some $C_1>0$ such that
\begin{align*}
qL_f\intomega\frac{\uep^{q-1}}{(2\delta-\wep)^\theta}\leq \intomega\frac{\uep^q}{(2\delta-\wep)^\theta}+C_1\intomega\frac{1}{(2\delta-\wep)^\theta}\leq\intomega\frac{\uep^q}{(2\delta-\wep)^\theta}+\frac{C_1}{\delta^\theta}\quad\text{for all }t>T.
\end{align*}
Plugging these two estimates into \eqref{eq:testing-uep-eq1} finally completes the proof.
\end{bew}

In order to get rid of the last gradient term appearing in the differential inequality of the previous lemma, we will rely on the decay property of $\wep$ to choose $\delta$ sufficiently small. The choice of $\delta$ is quite similar to the reasoning found in \cite[Lemma 3.5]{Lan-Longterm_M3AS16}, however, due to the more general form of $f$, here we will require some additional information to make successful use of the resulting differential inequality. 

\begin{lemma}
\label{lem:ev-lp-bound-uep}
Assume $\mu>0$, $\alpha>1+\sqrt{2}$ and $\beta>1$ and suppose that $\rho:=\min\{\alpha,\beta\}$ satisfies $\rho>\frac{\alpha+1}{\alpha-1}$ and that $r\in\CSp{1}{\bomega\times[0,\infty)}\cap\LSp{\infty}{\Omega\times(0,\infty)}\cap\LSp{1}{(0,\infty);\Lo[\infty]}$ fulfills \eqref{rprop2}. Then for all $q>1$ there exist $T>0$ and $C>0$ such that for each $\epsi\in(0,1)$ and all $t>T$ the solution $(\uep,\vep,\wep)$ of \eqref{approxprob} fulfills
\begin{align*}
\intomega\uep^q(\cdot,t)\leq C
\end{align*}
\end{lemma}

\begin{bew}
Given $q>1$ we fix a small $\theta>0$ satisfying
\begin{align}\label{eq:ev-lp-bound-uep-eq1}
4\big(q+2q(q-1)+q(q-1)^2\big)\theta<2(q-1)
\end{align}
and then pick $0<\delta<\min\big\{1,\theta,\frac{1}{4q}\big\}$. Since $\delta<\frac{1}{4q}$ implies $1<2(\theta+1-2q\delta)$ and $\delta<\theta$ entails $\frac{\delta}{\theta}<1$, we conclude from \eqref{eq:ev-lp-bound-uep-eq1} that
\begin{align*}
4\big(q+2q(q-1)\frac{\delta}{\theta}+q(q-1)^2\frac{\delta^2}{\theta^2}\big)\theta<4(q-1)(\theta+1-2q\delta),
\end{align*}
which, upon simple rearrangement, may also be expressed as
\begin{align}\label{eq:ev-lp-bound-uep-eq2}
\frac{\big(2q\theta+2q(q-1)\delta\big)^2}{4\big(\theta(\theta+1)-2q\theta\delta\big)}<q(q-1).
\end{align}
With these parameters fixed, we then employ Lemma \ref{lem:dec-wep-Linfty} to find $T>0$ such that
\begin{align*}
\big\|\wep(\cdot,t)\big\|_{\Lo[\infty]}<\delta\quad\text{for all }\epsi\in(0,1)\text{ and }t>T
\end{align*}
and conclude from Lemma \ref{lem:testing-uep-wep} and \eqref{eq:ev-lp-bound-uep-eq2} that there is some $C_1>0$ such that
\begin{align*}
\frac{\intd}{\intd t}\intomega\frac{\uep^q}{(2\delta-\wep)^\theta}+qK_f\intomega\frac{\uep^{q+\alpha-1}}{(2\delta-\wep)^\theta}\leq C_1\intomega\frac{\uep^q}{(2\delta-\wep)^\theta}+C_1\quad\text{for all }\epsi\in(0,1)\text{ and }t>T.
\end{align*}
Here, Hölder's inequality and Young's inequality entail the existence of $C_2>0$ such that
\begin{align*}
2C_1\intomega\frac{\uep^q}{(2\delta-\wep)^\theta}\leq 2C_1\bigg(\frac{|\Omega|}{\delta^\theta}\bigg)^{\frac{\alpha-1}{q+\alpha-1}}\bigg(\intomega\frac{\uep^{q+\alpha-1}}{(2\delta-\wep)^\theta}\bigg)^\frac{q}{q+\alpha-1}\leq qK_f\intomega\frac{\uep^{q+\alpha-1}}{(2\delta-\wep)^\theta}+C_2
\end{align*}
for all $\epsi\in(0,1)$ and $t>T$, so that upon letting $y_\epsi(t):=\intomega\frac{\uep^q(\cdot,t)}{(2\delta-\wep(\cdot,t))^\theta}$, $t>T$, we find that $y_\epsi$ satisfies
\begin{align}\label{eq:ev-lp-bound-uep-eq3}
y'_\epsi +C_1 y_\epsi\leq C_1+C_2=:C_3\quad\text{for all }\epsi\in(0,1)\text{ and }t>T.
\end{align}
To exploit the information contained in this differential inequality successfully, we note that according to Corollary \ref{cor:st-bounds-uep-l2} there is some $C_4>0$ such that 
$\int_T^{T+1}\!\!\intomega|\nabla\uep|^2\leq C_4$ and that hence for each $\epsi\in(0,1)$ we can find $t_\epsi\in(T,T+1)$ such that $\intomega|\nabla\uep(\cdot,t_\epsi)|^2\leq C_4.$ In view of the embedding $\W[1,2]\hookrightarrow\Lo[q]$ this also readily implies that there is some $C_5>0$ such that
\begin{align*}
\intomega\uep^q(\cdot,t_\epsi)\leq C_5\quad\text{for all }\epsi\in(0,1).
\end{align*}
From a combination with the ODE comparison argument from Lemma \ref{lem:diffineq-lemma} applied to \eqref{eq:ev-lp-bound-uep-eq3} we can therefore conclude the existence of $C_6>0$ such that
\begin{align*}
y_\epsi(t)\leq y_\epsi(t_\epsi)+\frac{C_3}{1-e^{-C_1}}\leq C_6\quad\text{for all }\epsi\in(0,1)\text{ and all }t>T+2>t_{\epsi}+1.
\end{align*}
Due to $(2\delta-\wep)^\theta\leq(2\delta)^\theta$ on $\Omega\times(T,\infty)$ we can finally establish that for all $\epsi\in(0,1)$
\begin{align*}
\intomega\uep^q=\intomega(2\delta-\wep)^\theta\frac{\uep^q}{(2\delta-\wep)^\theta}\leq(2\delta)^\theta\intomega\frac{\uep^q}{(2\delta-\wep)^\theta}= (2\delta)^\theta y_\epsi(t)\leq (2\delta)^\theta C_6\quad\text{on }(T+2,\infty),
\end{align*}
completing the proof.
\end{bew}

\setcounter{equation}{0} 
{\hypersetup{linkcolor=black}
\section{Eventual smoothness - Proof of Theorem \ref{theo:2}}\label{sec11:ev-smooth}}

With the eventual $L^q$ bounds of $\uep$ for any $q>1$ from the previous section at hand, we can now return to our conditional estimates of Section \ref{sec8:ev-cond-sob-reg} to lift the regularity information on the solution in an iterative procedure even higher, finally concluding in an eventual Hölder estimates for both $\uep$ and $\vep$. For this, however, we have to ensure that the spatio-temporal bound on $\vep$ in Lemma \ref{lem:spat-temp-alpha-beta} is sufficiently strong to work as a starting point for this procedure, which is the essential reason why we will have to restrict ourselves to values of $\beta>1+\sqrt{2}$. 

\begin{lemma}
\label{lem:hoelder-uep-wep}
Assume $\mu>0$, $\alpha>1+\sqrt{2}$ and $\beta>1+\sqrt{2}$ and that $r\in\CSp{1}{\bomega\times[0,\infty)}\cap\LSp{\infty}{\Omega\times(0,\infty)}\cap\LSp{1}{(0,\infty);\Lo[\infty]}$ fulfills \eqref{rprop2}. Then for any $p\in(1,\infty)$ there are $T>0$ and $C>0$ such that for each $\epsi\in(0,1)$ and all $t>T$ the solution $(\uep,\vep,\wep)$ of \eqref{approxprob} satisfies
\begin{align}\label{eq:sob-reg-for-all-p}
\intomega\vep^p(\cdot,t)+\int_t^{t+1}\!\!\big\|\uep\big\|_{\W[2,p]}^p+\int_t^{t+1}\!\!\big\|\wep\big\|_{\W[2,p]}^p\leq C.
\end{align}
Moreover, there exist $\gamma\in(0,1)$ and $C'>0$ such that for any $\epsi\in(0,1)$ and all $t>T$ the solution $(\uep,\vep,\wep)$ of \eqref{approxprob} satisfies
\begin{align}\label{eq:hoelder-reg-uep-wep}
\big\|\uep\big\|_{\CSpnl{1+\gamma,\frac{\gamma}{2}}{\bomega\times[t,t+1]}}+\big\|\wep\big\|_{\CSpnl{1+\gamma,\frac{\gamma}{2}}{\bomega\times[t,t+1]}}\leq C'.
\end{align}
\end{lemma}

\begin{bew}
The proof of \eqref{eq:sob-reg-for-all-p} is based on an inductive argument relying on the conditional eventual regularity estimates from Lemma \ref{lem:ev-max-sob}. We set $q_{(-1)}:=\beta-1$ and for $k\in\N_0$ define 
$q_{(k)}:=k+\beta$. We claim that for all $k\in\N$ we have $q_{(k-1)}>\max\{\beta-1,\frac{1}{\beta-2}\}$ and that there are $T_{(k)}>0$ and $C_{(k)}>0$ such that for each $\epsi\in(0,1)$ and all $t>T_{(k)}$
\begin{align}\label{eq:induction-eq1}
\int_t^{t+1}\!\big\|\uep\big\|_{\W[2,q_{(k-2)}]}^{q_{(k-2)}}+\int_t^{t+1}\!\big\|\wep\big\|_{\W[2,q_{(k)}]}^{q_{(k)}}+\intomega\vep^{1-\beta+q_{(k)}}(\cdot,t)+\int_t^{t+1}\!\!\intomega\vep^{q_{(k)}}\leq C_{(k)}.
\end{align}
It can easily be checked, that $\beta>1+\sqrt{2}$ implies $q_{(0)}=\beta>\max\{\beta-1,\frac{1}{\beta-2}\}$. Moreover, noting that $q_{(0)}-\beta+2=2<\beta$, we conclude from Lemma \ref{lem:spat-temp-alpha-beta}, Corollary \ref{cor:cond-spat-temp-est-W2q-wep} and Lemma \ref{lem:ev-lp-bound-uep} that the conditions of Lemma \ref{lem:ev-max-sob} are met for $q=q_{(0)}=\beta$ and that, hence, there are $T_{(1)}>0$ and $C_{(1)}>0$ such that
\begin{align*}
&\int_t^{t+1}\!\big\|\uep\big\|_{\W[2,q_{(-1)}]}^{q_{(-1)}}+\int_t^{t+1}\!\big\|\wep\big\|_{\W[2,q_{(1)}]}^{q_{(1)}}+\intomega\vep^{1-\beta+q_{(1)}}+\int_t^{t+1}\!\!\intomega\vep^{q_{(1)}}\\
=\ &\int_t^{t+1}\!\big\|\uep\big\|_{\W[2,\beta-1]}^{\beta-1}+\int_t^{t+1}\!\big\|\wep\big\|_{\W[2,\beta+1]}^{\beta+1}+\intomega\vep^{2}+\int_t^{t+1}\!\!\intomega\vep^{\beta+1}\leq C_{(1)}
\end{align*}
for each $\epsi\in(0,1)$ and all $t>T_{(1)}$, proving that the claim holds for $k=1$. Now, assuming that for an arbitrary but fixed $k\in\N$ we have $q_{(k-1)}>\max\{\beta-1,\frac{1}{\beta-2}\}$ and $T_{(k)}>0$, $C_{(k)}>0$ such that \eqref{eq:induction-eq1} holds for each $\epsi\in(0,1)$ and all $t>T_{(k)}$, we proceed to show the inductive step from $k$ to $k+1$. From the definition of $q_{(k-1)}$ we find that $q_{(k)}=q_{(k-1)}+1>\max\{\beta-1,\frac{1}{\beta-2}\}$ and by the induction hypothesis there is some $C'>0$ such that
\begin{align*}
\int_t^{t+1}\!\big\|\wep\big\|_{\W[2,\beta+k]}^{\beta+k}+\int_t^{t+1}\!\intomega\vep^{2+k}\leq\int_t^{t+1}\!\big\|\wep\big\|_{\W[2,\beta+k]}^{\beta+k}+\int_t^{t+1}\!\intomega\vep^{\beta+k}+C'\leq C_{(k)}+C'
\end{align*}
holds for each $\epsi\in(0,1)$ and all $t>T_{(k)}$. Together with Lemma \ref{lem:ev-lp-bound-uep} these estimates ensure that Lemma \ref{lem:ev-max-sob} becomes applicable for $q=q_{(k)}=\beta+k$. Hence, we can conclude that there are $T_{(k+1)}>0$ and $C_{(k+1)}>0$ such that \eqref{eq:induction-eq1} remains true when inserting $k+1$, implying that, actually, for any $k\in\N$ we can find $T_{(k)}>0$ and $C_{(k)}>0$ such that \eqref{eq:induction-eq1} holds for each $\epsi\in(0,1)$ and all $t>T_{(k)}$.
This proves the first assertion of the lemma due to $q_{(k)}>k+2$ for all $k\in\N$. From this we can also readily derive \eqref{eq:hoelder-reg-uep-wep} in view of known embedding results (e.g. \cite[Theorem 1.1]{Amann00}) by taking $k$ sufficiently large. 
\end{bew}

Since the previous lemma contains sufficiently strong information on the derivatives of $\uep$, we can now employ the maximal Sobolev estimate one final time to the second equation of \eqref{approxprob} in order to obtain an eventual Hölder estimate for $\vep$.

\begin{lemma}
\label{lem:hoelder-vep}
Assume $\mu>0$, $\alpha>1+\sqrt{2}$ and $\beta>1+\sqrt{2}$ and that $r\in\CSp{1}{\bomega\times[0,\infty)}\cap\LSp{\infty}{\Omega\times(0,\infty)}\cap\LSp{1}{(0,\infty);\Lo[\infty]}$ fulfills \eqref{rprop2}. Then there exist $\gamma\in(0,1)$, $T>0$ and $C>0$ such that for each $\epsi\in(0,1)$ and all $t>T$ the solution $(\uep,\vep,\wep)$ of \eqref{approxprob} satisfies
\begin{align*}
\big\|\vep\big\|_{\CSpnl{1+\gamma,\frac{\gamma}{2}}{\bomega\times[t,t+1]}}\leq C.
\end{align*}
\end{lemma}

\begin{bew}
Given $\theta>1$ we pick $\sigma\in(1,\theta)$ and note that according to Lemma \ref{lem:hoelder-uep-wep} we can fix $t_0>0$, and $C_1>0$ such that for all $t>t_0$ and each $\epsi\in(0,1)$ we have
\begin{align}\label{eq:hoelder-vep-eq1}
\big\|\nabla\uep\big\|_{\LSp{\infty}{\Omega\times(t,t+2)}}+\int_{t}^{t+2}\!\big\|\Delta\uep\big\|_{\Lo[2\theta]}^\theta\leq C_1
\end{align}
and $C_2>0$ such that
\begin{align}\label{eq:hoelder-vep-eq2}
\|\vep\|_{\LSp{\infty}{(t_0,\infty);\Lo[2\theta]}}+\|\vep\|_{\LSp{\infty}{(t_0,\infty);\Lo[\beta\theta]}}+\|\vep\|_{\LSp{\infty}{(t_0,\infty);\Lo[\theta]}}+\|\vep\|_{\LSp{\infty}{(t_0,\infty);\Lo[\sigma]}}\leq C_2
\end{align}
is valid for each $\epsi\in(0,1)$. Now, in order to prepare maximal Sobolev regularity arguments in a fashion similar to Lemma \ref{lem:ev-max-sob}, we let $\xi:=\xi_{t_0}$ again be given by Definition \ref{def:cutoff} and observe that $\xi\vep$ satisfies
\begin{align*}
\big(\xi\vep\big)_t=\Delta\big(\xi\vep\big)-\nabla\big(\xi\vep\big)\nabla\uep-\xi\vep\Delta\uep+\xi g(\vep)+\xi'\vep\quad\text{in }\Omega\times(t_0,\infty)
\end{align*}
with $\nabla(\xi\vep)\cdot\nu=0$ on $\romega\times(t_0,\infty)$ and $\xi\vep(\cdot,t_0)=0$ in $\Omega$. Hence, maximal Sobolev regularity estimates for the Neumann heat semigroup (\cite{HieberPruss-CPDE92,GigSohr91}) entail the existence of $C_3>0$ such that for each $\epsi\in(0,1)$ we may estimate
\begin{align*}
\int_{t_0}^{t_0+2}\!\big\|(\xi\vep)_t\big\|_{\Lo[\theta]}^\theta+\int_{t_0}^{t_0+2}\!\big\|\xi\vep\big\|_{\W[2,\theta]}^\theta
&\leq C_3\int_{t_0}^{t_0+2}\!\big\|\nabla(\xi\vep)\cdot\nabla\uep\big\|_{\Lo[\theta]}^\theta+C_3\int_{t_0}^{t_0+2}\!\big\|\xi\vep\Delta\uep\big\|_{\Lo[\theta]}^\theta\\
&\qquad+C_3\int_{t_0}^{t_0+2}\!\big\|\xi g(\vep)\big\|_{\Lo[\theta]}^\theta+C_3\int_{t_0}^{t_0+2}\!\big\|\xi'\vep\big\|_{\Lo[\theta]}^\theta.
\end{align*}
Making use of Hölder's inequality, \eqref{fgprop} and the fact that $\xi\leq1$ on $\R$, we find $C_4>0$ such that for each $\epsi\in(0,1)$
\begin{align}\label{eq:hoelder-vep-eq3}
&\int_{t_0}^{t_0+2}\!\big\|\xi\vep\big\|_{\W[2,\theta]}^\theta\nonumber\\
\leq\ &C_3\big\|\nabla\uep\big\|_{\LSp{\infty}{(\Omega\times(t_0,t_0+2)}}\!\int_{t_0}^{t_0+2}\!\!\intomega\big|\nabla(\xi\vep)\big|^\theta+C_3\!\int_{t_0}^{t_0+2}\!\big\|\vep\big\|_{\Lo[2\theta]}^\theta\big\|\Delta\uep\big\|_{\Lo[2\theta]}^\theta\\
&\quad+ C_4\int_{t_0}^{t_0+2}\!\!\intomega\vep^{\beta\theta}+ C_4+ C_3\|\xi'\|_{\LSp{\infty}{\R}}^\theta\int_{t_0}^{t_0+2}\!\!\intomega\vep^\theta.\nonumber
\end{align}
Since a combination of the \GNI\ and Young's inequality provides $C_5,C_6>0$ and $a=\frac{\frac{1}{2}+\frac{1}{\sigma}-\frac{1}{\theta}}{1+\frac{1}{\sigma}-\frac{1}{\theta}}<1$ such that
\begin{align*}
C_1C_3\|\phi\|_{\W[1,\theta]}^\theta\leq C_5\|\phi\|_{\W[2,\theta]}^{a\theta}\|\phi\|_{\Lo[\sigma]}^{(1-a)\theta}\leq \frac{1}{2}\|\phi\|_{\W[2,\theta]}^\theta+C_6\|\phi\|_{\Lo[\sigma]}^\theta\quad\text{for all }\phi\in\W[2,\theta],
\end{align*}
we conclude from \eqref{eq:hoelder-vep-eq1}, \eqref{eq:hoelder-vep-eq2} and \eqref{eq:hoelder-vep-eq3} that for each $\epsi\in(0,1)$ we have
\begin{align*}
\frac{1}{2}\int_{t_0}^{t_0+2}\!\big\|\xi\vep\big\|_{\W[2,\theta]}^\theta\leq 2C_2^\theta C_6+ C_1 C_2^\theta C_3+2C_2^{\beta\theta} C_4+C_4+2C_2^{\theta}C_3\|\xi'\|_{\LSp{\infty}{\R}}^\theta.
\end{align*}
Hence, noting that $\xi\equiv 1$ on $(t_0+1,t_0+2)$, we can for any $\theta>1$ find some $C_7>0$ such that for each $\epsi\in(0,1)$ and all $t>T:=t_0+1$ we have
\begin{align*}
\int_t^{t+1}\!\|\vep\|_{\W[2,\theta]}^\theta\leq C_7.
\end{align*}
The asserted Hölder regularity now follows from an application of the known embedding results e.g. presented in \cite[Theorem 1.1]{Amann00}.
\end{bew}

With the prepared regularity information we can now rely on the Arzelà--Ascoli theorem and parabolic regularity theory to verify that the global generalized solution obtained in Theorem \ref{theo:1} actually solves \eqref{forager} classically at least after some smoothing time $\Td>0$. Reasoning along these lines has previously been utilized in similar settings and can e.g. be found in \cite[Lemma 7.13]{win_chemonavstokesfinal_TransAm17}, \cite[Lemma 3.12]{Lan-Longterm_M3AS16} and \cite[Lemma 7.5]{TB2019_stokeslimit}.

\begin{proof}[\textbf{Proof of Theorem \ref{theo:2}}:]
Making use of Lemma \ref{lem:hoelder-uep-wep} and Lemma \ref{lem:hoelder-vep} we can pick some $\gamma'\in(0,1)$, $T>0$ and $C_1>0$ such that for each $\epsi\in(0,1)$ and any $t>T$
\begin{align*}
\big\|\uep\big\|_{\CSp{1+\gamma',\frac{\gamma'}{2}}{\bomega\times[t,t+1]}}+
\big\|\vep\big\|_{\CSp{1+\gamma',\frac{\gamma'}{2}}{\bomega\times[t,t+1]}}+
\big\|\wep\big\|_{\CSp{1+\gamma',\frac{\gamma'}{2}}{\bomega\times[t,t+1]}}\leq C_1.
\end{align*}
Denoting by $(\epsi_j)_{j\in\N}$ the sequence obtained in Proposition \ref{prop:conv} and with $(u,v,w)$ the corresponding limit object, we can draw on the Arzelà--Ascoli theorem to find that for some $\gamma_0\in(0,\gamma')$
\begin{align*}
\uep\to u,\quad \vep\to v,\quad \wep\to w\quad\text{in }\CSp{1+\gamma_0,\frac{\gamma_0}{2}}{\bomega\times[t,t+1]}
\end{align*}
along a subsequence of $(\epsi_j)_{j\in\N}$. In particular, $u$ and $v$ are bounded on $\bomega\times[t,t+1]$, which in light of $f,g\in\CSp{1}{[0,\infty)}$, implies that
\begin{align*}
f(u),\  g(v)\in\CSp{\gamma',\frac{\gamma_0}{2}}{\bomega\times[t,t+1]}.
\end{align*}
Letting $\xi_1:=\xi_T$ be provided by Definition \ref{def:cutoff} we find that $W=\xi_1 w$ solves the problem
\begin{align*}
W_t=\Delta W + H,\quad W(T)=0,\quad\frac{\partial W}{\partial \nu}\big|_{\romega}=0,
\end{align*}
with $H=-(u+v)\xi_1 w-\mu\xi_1 w+r\xi_1+\xi_1' w\in\CSp{\gamma_0,\frac{\gamma_0}{2}}{\bomega\times(T,\infty)}$. Hence, parabolic theory (\cite[III.5.1 and IV.5.3]{lsu}) tells us that there is some $\gamma_1\in(0,\gamma_0)$ such that 
\begin{align}\label{eq:theo2-eq1}
w\in\CSp{2+\gamma_1,1+\frac{\gamma_1}{2}}{\bomega\times[T+1,\infty)}.
\end{align}
Similarly, setting $\xi_2:=\xi_{T+1}$ we notice that $U=\xi_2 u$ is a solution of the problem
\begin{align*}
U_t=\Delta U- A\cdot \nabla U+ B,\quad U(T+1)=0,\quad\frac{\partial U}{\partial\nu}\big|_{\romega}=0,
\end{align*}
where $A=\nabla w\in\CSp{\gamma_1,\frac{\gamma_1}{2}}{\bomega\times(T+1,\infty)}$ and $B=-\xi_2 u\Delta w+\xi_2 f(u)+\xi_2'u\in\CSp{\gamma_1,\frac{\gamma_1}{2}}{\bomega\times(T+1,\infty)}$ again imply by parabolic theory (\cite[III.5.1 and IV.5.3]{lsu}) that there is $\gamma_2\in(0,\gamma_1)$ such that 
\begin{align}\label{eq:theo2-eq2}
u\in\CSp{2+\gamma_2,1+\frac{\gamma_2}{2}}{\bomega\times[T+2,\infty)}.
\end{align}
Repeating the same arguments a third time with $\xi_3:=\xi_{T+2}$ for the solution $V=\xi_3 v$ of
\begin{align*}
V_t=\Delta V-A'\cdot\nabla V+B',\quad V(T+2)=0,\quad\frac{\partial V}{\partial \nu}\big|_{\romega}=0,
\end{align*}
with $A'=\nabla u\in\CSp{\gamma_2,\frac{\gamma_2}{2}}{\bomega\times(T+3,\infty)}$ and $B'=-\xi_3 v\Delta u+\xi_3 g(v)+\xi_3' v\in\CSp{\gamma_2,\frac{\gamma_2}{2}}{\bomega\times(T+2,\infty)}$ lastly provides $\gamma_3\in(0,\gamma_2)$ such that 
\begin{align}\label{eq:theo2-eq3}
v\in\CSp{2+\gamma_3,1+\frac{\gamma_3}{2}}{\bomega\times[T+3,\infty)}.
\end{align}
As a consequence of \eqref{eq:theo2-eq1}, \eqref{eq:theo2-eq2}, \eqref{eq:theo2-eq3} we readily obtain the regularity assertions of the theorem, which, upon letting $\Td:=T+3$, imply that $(u,v,w)$ solves \eqref{forager} classically in $\Omega\times(\Td,\infty)$.
\end{proof}

\section*{Acknowledgements}
The author acknowledges support of the {\em Deutsche Forschungsgemeinschaft} in the context of the project
  {\em Emergence of structures and advantages in cross-diffusion systems (project no. 411007140)}.

\footnotesize{
\setlength{\bibsep}{2pt plus 0.5ex}

}

\end{document}